\documentclass[a4paper,onecolumn,oneside,11pt,leqno]{article}

\usepackage{amsmath, amsfonts,amssymb}
\usepackage{theorem, euscript,array,enumerate}
\usepackage{amsfonts,mathrsfs, latexsym}
\usepackage[dvips]{graphicx}
\usepackage[dvips]{graphics}
\usepackage[dvips]{epsfig}

\theoremheaderfont{\bfseries} \theoremstyle{plain}
\theorembodyfont{\upshape}
\newcommand{\esp}{\mathbb{E}} %esperance

\newcommand{\indic}{1\!\!1} %indicatrice
\newcommand{\proba}{\mathbb{P}}

\newcommand{\var}{\mathrm{Var}}

\newcommand{\R}{\mathbb{R}}

\theoremstyle{plain}

\theorembodyfont{\upshape\it}
\newtheorem{Thm}{\bf Theorem}[section]
\newtheorem{Pro}[Thm]{\bf Proposition}
\newtheorem{Lem}[Thm]{\bf Lemma}
\newtheorem{Cor}[Thm]{\bf Corollary}
\theorembodyfont{\upshape}

\newtheorem{Def}[Thm]{Definition}

\newtheorem{Rem}[Thm]{Remark}

\newenvironment{proof}[1][{\it Proof.}]{\begin{trivlist}
\item[\hskip \labelsep {\bfseries #1}]}
{ {\quad}\hfill $\Box$\end{trivlist}\vskip
-0.2 cm}

\begin{document}
\thispagestyle{empty}
\title{Zero bias transformation and asymptotic expansions}
\author{Ying Jiao\thanks{Laboratoire de probabilit\'es et mod\`eles
al\'eatoires, Universit\'e Paris 7,
jiao@math.jussieu.fr. 
} }
\date{\today}
\maketitle
\begin{abstract}
We apply the zero bias transformation to deduce a recursive asymptotic expansion formula for expectation of functions of sum of independent random variables  in terms of normal expectations and we discuss the remainder term estimations. 
%The main tools are a so-called reverse Taylor formula and a concentration inequality.
\vspace{3mm}

\end{abstract} 

\section{Introduction}\label{Sec:introduction}
\hskip\parindent
Zero bias transformation has been introduced by Goldstein and Reinert \cite{GR1997}
in the framework of Stein's method. By the fundamental works of Stein \cite{St1972, St1986}, we know that
a random variable (r.v.) $Z$ with mean zero follows the normal distribution $N(0,\sigma^2)$
if and only if $\esp[Zf(Z)]=\sigma^2\esp[f'(Z)]$ for any Borel function $f$ such that both sides of the
equality are well defined. More generally,
for any r.v. $X$ with mean zero and finite variance $\sigma^2>0$, a r.v.
$X^*$ is said to have the zero biased distribution
of $X$ if the equality
\begin{equation}\label{Equ:zero bais distr}\esp[Xf(X)]=\sigma^2\esp[f'(X^*)]\end{equation}
holds for any differentiable function $f$ such that \eqref{Equ:zero bais distr} is well defined.
 So combined with the {\it Stein's equation}
$xf(x)-\sigma^2f'(x)=h(x)-\Phi_\sigma(h)$ where $h$ is a given function and
$\Phi_\sigma(h)$ denotes the expectation of $h$
under $N(0,\sigma^2)$, we have
$$\esp[h(X)]-\Phi_\sigma(h)=\esp[Xf_h(X)-\sigma^2f_h'(X)]=\sigma^2\esp[f_h'(X^*)-f'_h(X)]$$where
$f_h$ is the solution of Stein's equation given by
\begin{equation}\label{Equ:solution of Stein's equation}
f_h:=\frac{1}{\sigma^2\phi_\sigma(x)}\int_x^\infty
(h(t)-\Phi_\sigma(h))
\phi_\sigma(t)\,dt
%=-\frac{1}{\sigma^2\phi_\sigma(x)} \int^x_{-\infty} \overline h(t)
%\phi_\sigma(t)\,dt
.\end{equation}
An important remark is that
$X^*$ need not be independent of $X$ (\cite{GR1997}, see also \cite{Go2007}).
In fact, let $W=X_1+\cdots+X_n$ be sum of independent mean zero
random variables,
Goldstein and Reinert proposed the construction $W^*:=W^{(I)}+X_I^*$
where, for any $i\in\{1,\cdots,n\}$, $W^{(i)}:=W-X_i$ and $X_i^*$ is
independent of $W^{(i)}$, and $I$
is a random index valued in $\{1,\cdots,n\}$ which is
independent of $(X_1,\cdots,X_n,X_1^*,\cdots,X_n^*)$
and satisfies $\mathbb P(I=i)=\sigma_i^2/\sigma_W^2$
with $\sigma_i^2$ being the variance of $X_i$ and $\sigma_W^2$
that of $W$.
We observe that the above construction
of zero bias transformation is quite similar to
Lindeberg method except that, in zero bias
transformation, we consider an average of punctual
substitutions of $X_i$ by $X_i^*$; while in Lindeberg
method, we substitute progressively $X_i$ by central
normal distribution with the same variance.

The asymptotic expansion of
expectations of the form $\esp[h(W)]$ is a classical
topic in central limit theorems. Using Stein's method,
Barbour \cite{Ba1986,Ba1987} has obtained a full
expansion of $\esp[h(W)]$ for sufficiently regular
function $h$. Compared to the classical Edgeworth
expansion (see \cite[ ChapV]{Pe1975}, also
\cite{Ro2005}), the results of \cite{Ba1986} do not
require the distribution of $X_i$ to be smooth;
however, as a price paid, we  need some suitable
regularity conditions on the function $h$. The result of
\cite{Ba1986} can also be compared to those in
\cite{Hi1977,GH1978} using Fourier transform. The key
point of Barbour's method is a Taylor type formula with
cumulant coefficients, which allows to write the
difference  $\mathbb
E[Wf(W)]-\sigma_W^2\mathbb E[f'(W)]$ as a series which
involves cumulants of order $\ge 3$ and
to iterate the procedure of replacing $W$-expectations by
normal expectations until the desired order. It has been pointed out in \cite{Ro2005} that
the key formula of Barbour can also be obtained by Fourier transform.

Zero bias transformation have been used in \cite{EJ2008} to obtain a first order
correction term for the normal approximation of $\esp[h(W)]$,
where the motivation was to find a rapid numerical method for large-sized credit
derivatives. The function of interest is the so-called
{\it call function} in finance: $h(x)=(x-k)_+$ where
$k$ is a real number. Since such $h$ is only absolutely
continuous, the function $f_h$ is not regular enough to
have the third order derivative. To
achieve the estimation, the authors have used a
conditional expectation technique, together with a
concentration inequality due to Chen and Shao
\cite{CS2001,CS2005}.

The main difficulty in generalizing the result in \cite{EJ2008}
to obtain a full expansion of $\esp[h(W)]$ is that
$W$ and $W^*-W$ are not independent. In fact, if we consider
the Taylor expansion of $f_h'(W^*)$ at $W$ and then
apply the expectation, there appear terms of the form
$\esp[f_h^{(l)}(W)(W^*-W)^k]$, where $f^{(l)}$ denotes the $l^{\text{th}}$-order derivative of  $f$. For the first order
expansion in \cite{EJ2008}, the conditional expectation argument allows
us to replace $\esp[f''_h(W)(W^*-W)]$ by
$\esp[f''_h(W)]\esp[W^*-W]$ and put the
covariance in the error term. However, in higher order
expansion, the error term could no longer contain such
covariances.
An alternative way is to consider the
Taylor expansion of
$\esp[f_h'(W^*)-f_h'(W)]$ at
$W^{(i)}$. As $X_i^*$ is independent of $W^{(i)}$,
there is no crossing term. However, the expectations of
the form $\esp[f_h^{(l)}(W^{(i)})]$ appear, which make
it difficult to apply the recurrence procedure. To overcome this
difficulty, we propose a so-called reverse Taylor formula which enables us to
replace $\esp[f_h^{(l)}(W^{(i)})]$ by expectation of
functions of $W$, up to an error term.

%=======================================
%\subsection{Reverse Taylor formula }
%=======================================
%\hskip\parindent

Let $N$ be a positive integer, $X$ and $Y$ be two
independent random variables such that $Y$ has up to
$N^{\mathrm{th}}$ order moments, and $f$ be an
$N^{\mathrm{th}}$ order differentiable function such
that $f^{(k)}(X)$ and $f^{(k)}(X+Y)$ are integrable for
any $k=0,\cdots, N$.
We define the notation $m_Y^{(k)}:=\esp[Y^k]/k!$. Denote by
$\delta_N(f,X,Y)$ the error term in the
$N^{\mathrm{th}}$ order Taylor expansion of
$\esp[f(X+Y)]$. Namely,
\begin{equation}\delta_N(f,X,Y):=
\esp[f(X+Y)]-\sum_{k=0}^Nm_Y^{(k)}\esp[f^{(k)}(X)].\end{equation}
Recall  that for any $N\geq 1$,
\begin{equation}\label{Equ:remaining term of Taylor expansion}
\delta_N(f,X,Y)=\frac{1}{(N-1)!}
\int_0^1(1-t)^{N-1}\esp\Big[\big(f^{(N)}(X+tY)-f^{(N)}(X)\big)Y^{N}\Big]\,dt
\end{equation}
provided that the term on the right side is well
defined. This is a consequence of the classical Taylor formula in its integral form (e.g. \cite{La2001}).

The so-called {reverse Taylor formula} gives an expansion of $\esp[f(X)]$
in terms of expectations of functions of $X+Y$ and of 
moments of $Y$. We would like to note that, in the expansion formula
\eqref{Equ:reverse taylor fomrul}, the variables $X+Y$
and $Y$ are not independent.
We specify some notation and
conventions. First of all, $\mathbb N_*:=\mathbb
N\setminus\{0\}$ denotes the set of strictly positive
integers. For any integer $d\ge 1$ and any
$\mathbf{J}=(j_l)_{l=1}^d\in\mathbb N_*^d$,
$|\mathbf{J}|$ is defined as $j_1+\cdots+j_d$, and
$m_Y^{(\mathbf{J})}:=m_Y^{(j_1)}\cdots m_Y^{(j_d)}$. By
convention, $\mathbb N_*^0 $ denotes the set
$\{\emptyset\}$ of the empty vector, $|\emptyset|=0$
and $m_Y^{(\emptyset)}=1$.

\begin{Pro}\label{Lem:reverse taylor formula}{\bf(Reverse Taylor formula)}
With the  above notation, the
equality
\begin{equation}\label{Equ:reverse taylor fomrul}
\esp[f(X)]=\sum_{d\ge
0}(-1)^d\hspace{-4mm}\sum_{\mathbf{J}\in\mathbb
N^d_*,\, |\mathbf{J}|\le N
}\hspace{-4mm}m_Y^{(\mathbf{J})}\esp[f^{(|\mathbf{J}|)}(X+Y)]+
\varepsilon_N(f,X,Y)
\end{equation}
holds, where $\varepsilon_N(f,X,Y)$ is defined as
\begin{equation}\label{equ:reverse taylor remainder}
\varepsilon_N(f,X,Y)=-\sum_{d\ge 0}(-1)^d
\hspace{-4mm}\sum_{\mathbf{J}\in\mathbb N^d_*,\,
|\mathbf{J}|\le N
}\hspace{-4mm}m_Y^{(\mathbf{J})}\delta_{N-|\mathbf{J}|}(f^{(|\mathbf{J}|)},X,Y).
\end{equation}
\end{Pro}

%======================
%\subsection{Formal expansion}
%======================

%\hskip\parindent
The main result of this paper is an expansion formula for the sum of independent random variables.
We present below its formal form without giving precise conditions on the function
and on the summand variables (this will be done in Section 3). The methodology appeals to the zero bias transformation.
From now on, we consider a family of independent random variables $X_i \;(i=1,\cdots,n)$,
with mean zero and finite variance $\sigma_i^2>0$. Let $W=X_1+\cdots+X_n$ and $\sigma_W^2=\var(W)$.
Denote by $X_i^*$ a random variable which follows the zero-biased distribution of $X_i$
and which is independent of $W^{(i)}:=W-X_i$.
%Let $I$ be a random index valued in $\{1,\cdots,n\}$
%independent of $(X_1,\cdots,X_n,X_1^*,\cdots,X_n^*)$ and such that $\proba(I=i)=\sigma_i^2/\sigma_W^2$ for any $i$.

\begin{Thm}\label{Thm:main theorem}
Assume that $X_1,\cdots,X_n$ and the function $h$ are sufficiently
good (in a sense that we shall
precise later). Then, for any integer $N\geq 0$, $\esp[h(W)]$ can be written as the
sum of two terms $C_N(h)$ and $e_N(h)$, with
$C_0(h)=\Phi_{\sigma_W}(h)$ and
$e_0(h)=\esp[h(W)]-\Phi_{\sigma_W}(h)$, and recursively for $N\geq 1$,
\begin{equation}\label{Equ:recursive CN}
C_N(h)=C_0(h)+\sum_{i=1}^n\sigma_i^2\sum_{d\ge
1}(-1)^{d-1}\hspace{-5mm}\sum_{\mathbf{J}\in\mathbb
N_*^d,\, |\mathbf{J}|\le N}\hspace{-4mm}
m_{X_i}^{(\mathbf J^\circ)}\big(m_{X_i^*}^{(\mathbf
J^\dagger)}-m_{X_i}^{(\mathbf{J}^\dagger)}\big)
C_{N-|\mathbf{J}|}(f_h^{(|\mathbf{J}|+1)}),
\end{equation}
\begin{equation}\label{Equ:recursive eN}\begin{split}
e_N(h)&=\sum_{i=1}^n\sigma_i^2\bigg[ \sum_{d\ge
1}(-1)^{d-1}\hspace{-5mm}\sum_{\mathbf{J}\in\mathbb
N_*^d,\,|\mathbf{J}|\le
N}\hspace{-4mm}m_{X_i}^{(\mathbf
J^\circ)}\big(m_{X_i^*}^{(\mathbf
J^\dagger)}-m_{X_i}^{(\mathbf{J}^\dagger)}\big)e_{N-|\mathbf{J}|}(f_h^{(|\mathbf{J}|+1)})\\&
\qquad
+\sum_{k=0}^N\varepsilon_{N-k}(f_h^{(k+1)},W^{(i)},X_i)m_{X_i^*}^{(k)}
+\delta_N(f_h',W^{(i)},X_i^*)\bigg],
\end{split}\end{equation}
where for any integer $d\ge 1$, and any
$\mathbf{J}\in\mathbb N_*^d$, $\mathbf{J}^\dagger\in
\mathbb N_*$ denotes the last coordinate of
$\mathbf{J}$, and $\mathbf{J}^\circ$ denotes the
element in $\mathbb N_*^{d-1}$ obtained from
$\mathbf{J}$ by omitting the last coordinate.
\end{Thm}

In view of the classical formula relating the cumulants
and moments, % (see \cite{Pe1975} page 9 (2.3)),
our principal term $C_N(h)$ is similar to that
obtained by Barbour. Note that in $C_N(h)$, there appear
normal expectations of iteration of operators which are of the form
$g\mapsto f_g^{(l)}$ acting on $h$.  As pointed
out by Barbour \cite[p.294]{Ba1986}, such
expectation can be expressed as expectation of $h$
multiplied by a Hermite polynomial.
 
The proof of the equality $\esp[h(W)]=C_N(h)+e_N(h)$
is based on the reverse
Taylor formula and the zero bias transformation. It is important to
precise the conditions under which all terms
in the formal expansion are well defined. 
Moreover, we also need to show that $e_N(h)$
is ``small' enough as an error term. In our results,
the error term $e_N(h)$ is expressed in a
recursive way so that it is actually a linear
combination of remainders of Taylor and reverse
Taylor formulas and can be thus estimated.
A key ingredient in the estimation is a concentration inequality which provides upper bound for
$\proba(a\leq W\leq b)$  involving exponent $\leq 1$ of the interval length $(b-a)$, i.e. $(b-a)^\alpha$ with $0<\alpha\leq 1$. This allows to us to obtain, under relatively mild moment conditions on $X_i$'s than those in 
\cite{Ba1986}, estimations for the Taylor and reverse Taylor remainders. For example, as a consequence of Theorem \ref{Thm:main theorem} and the remainder estimations, we recover a classical result, initially obtained by using Fourier transform, asserting that if $X_1,\cdots,X_n$ are i.i.d. random variables with mean zero and variance $\sigma^2>0$, which admit $(2+\alpha)^{\text{th}}$ order moments, then the law of $(X_1+\cdots+X_n)/\sqrt n$ converges to $N(0,\sigma^2)$ and that the convergence speed is of order $({1}/{\sqrt n)^\alpha}$.

The rest
of the paper is organized as follows.  We
firstly prove the reverse Taylor formula and the formal expansion in
Section 2. In Section 3, we introduce the admissible function space and discuss the
conditions on $h$ and on $X_i$'s; this is inspired by ideas in \cite{Ba1986} and we can in addition include
some more irregular functions.
We then restate the main expansion result in this context. Section 4 is
devoted to error estimations. 
Finally, some
technical proofs are left in Appendix.

%================================================
\section{Reverse Taylor formula and formal expansion}
%================================================

\hskip\parindent To prove Proposition \ref{Lem:reverse taylor formula},
the main point is to replace
$\esp[f^{(|\mathbf{J}|)}(X+Y)]$ by its classical Taylor
expansion of $(N-|\mathbf{J}|)^{\mathrm{th}}$ order, so that all summand terms are of the same order and
some of them can be cancelled off progressively.

\begin{proof}[Proof of Proposition \ref{Lem:reverse taylor formula}]
We replace $\esp[f^{(|\mathbf{J}|)}(X+Y)]$ on the right
side of \eqref{Equ:reverse taylor fomrul} by
\[\sum_{k=0}^{N-|\mathbf{J}|}m_Y^{(k)}\esp[f^{(|\mathbf{J}|+
k)}(X)]+
\delta_{N-|\mathbf{J}|}(f^{(|\mathbf{J}|)},X,Y)\] and
observe that the sum of terms containing $\delta$
vanishes with $\varepsilon_N(f,X,Y)$. Hence we obtain
that the right side of \eqref{Equ:reverse taylor
fomrul} equals
\[\sum_{d\ge
0}(-1)^d\hspace{-4mm}\sum_{\mathbf{J}\in\mathbb
N^d_*,\, |\mathbf{J}|\le N
}\hspace{-4mm}m_Y^{(\mathbf{J})}\sum_{k=0}^{N-|\mathbf{J}|}m_Y^{(k)}\esp[f^{(|\mathbf{J}|+k)}(X)]
\]
If we split the last sum for $k=0$ and for $1\le
k\le N-|\mathbf{J}|$ respectively, the formula above
can be written as
\begin{equation}\label{Equ:proof of retf}\sum_{d\ge
0}(-1)^d\hspace{-4mm}\sum_{\mathbf{J}\in\mathbb
N^d_*,\, |\mathbf{J}|\le N
}\hspace{-4mm}m_Y^{(\mathbf{J})}\esp[f^{(|\mathbf{J}|)}(X)]
+\sum_{d\ge
0}(-1)^d\hspace{-4mm}\sum_{\mathbf{J}\in\mathbb
N^d_*,\, |\mathbf{J}|\le N
}\hspace{-4mm}m_Y^{(\mathbf{J})}\sum_{k=1}^{N-|\mathbf{J}|}m_Y^{(k)}\esp[f^{(|
\mathbf{J}|+k)}(X)].
\end{equation}
We make the index changes
$\mathbf{J}'=(\mathbf{J},k)$ and $u=d+1$ in the second
part, we find that it is just
\[\sum_{u\ge 1}(-1)^{u-1}\hspace{-4mm}
\sum_{\mathbf{J'}\in\mathbb N^u_*,\, |\mathbf{J'}|\le N
}\hspace{-4mm}m_Y^{(\mathbf{J'})}\esp[f^{(|\mathbf{J'}|)}(X)].\]
Thus, the terms in the first and the second parts of \eqref{Equ:proof of retf}
cancel out except the one of index
$d=0$ in the first part,
which proves the proposition.
\end{proof}

Using Proposition \ref{Lem:reverse taylor formula},
we prove below the formal equality
$\esp[h(W)]=C_N(h)+e_N(h)$ by induction on $N$.

\begin{proof}[Proof of Theorem \ref{Thm:main theorem} (formal part)]
The equality $\esp[h(W)]=C_0(h)+e_0(h)$ holds by
definition. In the following, we assume that the
equality $\esp[h(W)]=C_k(h)+e_k(h)$ has been verified
for any $k\in\{0,\cdots,N-1\}$ and for any good enough
function $h$. By Stein's equation, $\esp[h(W)]-C_0(h)$
is equal to
\[\sigma_W^2\esp[f_h'(W^*)-f_h'(W)]=
\sum_{i=1}^n\sigma_i^2\Big(\esp[f_h'(W^{(i)}+X_i^*)]-\esp[f_h'(W)]\Big).\]
Consider the following Taylor expansion
\[\esp[f_h'(W^{(i)}+X_i^*)]=\sum_{k=0}^Nm_{X_i^*}^{(k)}
\esp[f_h^{(k+1)}(W^{(i)})]+\delta_N(f_h',W^{(i)},X_i^*).\]
By replacing $\esp[f_h^{(k+1)}(W^{(i)})]$ in the above
formula
 by its
$(N-k)^{\mathrm{th}}$ reverse Taylor expansion, we
obtain that $\esp[f_h'(W^{(i)}+X_i^*)]$ equals
\[\sum_{k=0}^Nm_{X_i^*}^{(k)}\bigg[\sum_{d\ge 0}
(-1)^d\hspace{-2mm}\sum_{\begin{subarray}{c}\mathbf{J}\in\mathbb
N_*^d\\|\mathbf{J}|\le
N-k\end{subarray}}\hspace{-2mm}m_{X_i}^{(\mathbf{J})}
\esp[f_h^{(|\mathbf{J}|+k+1)}(W)]+\varepsilon_{N-k}
(f_h^{(k+1)},W^{(i)},X_i)\bigg]+\delta_N(f_h',W^{(i)},X_i^*).\]
Note that the term with indexes $k=d=0$ in the sum inside the bracket is
$\esp[f_h'(W)]$. Therefore
$\esp[f_h'(W^{(i)}+X_i^*)]-\esp[f_h'(W)]$ can be
written as the sum of the following three parts
\begin{gather}\label{Equ:first part}
\sum_{k=1}^n m_{X_i^*}^{(k)}\sum_{d\ge 0}
(-1)^d\hspace{-4mm}\sum_{\mathbf{J}\in\mathbb
N_*^d,\,|\mathbf{J}|\le
N-k}\hspace{-4mm}m_{X_i}^{(\mathbf{J})}
\esp[f_h^{(|\mathbf{J}|+k+1)}(W)],\\
\label{Equ:second part}\sum_{d\ge 1}
(-1)^d\hspace{-4mm}\sum_{\mathbf{J}\in\mathbb
N_*^d,\,|\mathbf{J}|\le
N}\hspace{-4mm}m_{X_i}^{(\mathbf{J})}
\esp[f_h^{(|\mathbf{J}|+1)}(W)],\\
\label{Equ:third part}
\sum_{k=0}^Nm_{X_i^*}^{(k)}\varepsilon_{N-k}(f_h^{(k+1)},W^{(i)},X_i)
+\delta_N(f_h',W^{(i)},X_i^*).
\end{gather}
By interchanging summations and then taking the index
changes $\mathbf{K}=(\mathbf{J},k)$ and $u=d+1$, we
obtain
\[\eqref{Equ:first part}=\sum_{u\ge 1}
(-1)^{u-1}\hspace{-4mm}\sum_{\mathbf{K}\in\mathbb
N_*^u,\,|\mathbf{K}|\le
N}\hspace{-4mm}m_{X_i}^{(\mathbf{K^\circ})}m_{X_i^*}^{(\mathbf{K}^\dagger)}
\esp[f_h^{(|\mathbf{K}|+1)}(W)].
\] As the equality $m_{X_i}^{(\mathbf{J})}=m_{X_i}^{(\mathbf{J}^\circ)}
m_{X_i}^{\mathbf{J}^\dagger}$ holds for any
$\mathbf{J}$, \eqref{Equ:first part}+\eqref{Equ:second
part} simplifies as
\[\sum_{d\ge 1}
(-1)^{d-1}\hspace{-4mm}\sum_{\mathbf{J}\in\mathbb
N_*^d,\,|\mathbf{J}|\le
N}\hspace{-4mm}m_{X_i}^{(\mathbf{J^\circ})}
\Big(m_{X_i^*}^{(\mathbf{J}^\dagger)}-m_{X_i}^{(\mathbf{J}^\dagger)}\Big)
\esp[f_h^{(|\mathbf{J}|+1)}(W)].\] By the hypothesis of
induction, we have
$\esp[f_h^{(|\mathbf{J}|+1)}(W)]=C_{N-|\mathbf{J}|}(f_h^{(|\mathbf{J}|+1)})+e_{N-|\mathbf{J}|}(f_h^{(|\mathbf{J}|+1)})$,
so the equality $\esp[h(W)]=C_N(h)+e_N(h)$ follows from
\eqref{Equ:recursive CN} and \eqref{Equ:recursive eN}.
\end{proof}

%=========================================
\section{Admissible function space}
%=========================================
\hskip\parindent In this section, we describe the function set for which we can
make the $N^{\text{th}}$ order expansion in Theorem \ref{Thm:main theorem}.
We need conditions on regularity and on the increasing speed at infinity of the
function $h$. Actually, from \eqref{Equ:recursive CN} and \eqref{Equ:recursive eN},
we are concerned with the $(N-k)^{\text{th}}$ order expansion
of $f_h^{(k+1)}$ for $k=1,\cdots,N$. So it would be natural to expect that $f'_h$ still belongs to this set.
Then by a recursive procedure, all terms will be well defined.

Recall (\cite{KF1974}, Chapter VI) that any function
$g$ on $\mathbb R$ which is locally of finite variation
can be uniquely decomposed into the sum of a function
of pure jump and a continuous function locally of
finite variation and vanishing at the origin. That is, $g=g_c+g_d$
where $g_c$ is called the {\it continuous part} of $g$ and
$g_d$ is the {\it purely
discontinuous part}.

Let  $\alpha\in(0,1]$ and $p\ge 0$
be two real numbers. For any function $f$ on
$\mathbb R$, the following quantity has been defined by Barbour in
\cite{Ba1986}:
\begin{equation}\label{equ:f norm alpha p}\|f\|_{\alpha,p}:=\sup_{x\neq y}\frac{|f(x)-f(y)|}{|x-y|^\alpha(
1+|x|^p+|y|^p)} .\end{equation} The finiteness of
this quantity implies that the function $f$ is locally 
$\alpha$-Lipschitz, and the increasing speed of $f$ at infinity is
at most of order $|x|^{\alpha+p}$. All
functions $f$ such that $\|f\|_{\alpha,p}<+\infty$ forms a
vector space over $\R$, and $\|\cdot\|_{\alpha,p}$ is a norm on
it. We list below several properties of
$\|\cdot\|_{\alpha,p}$, which will be useful
afterwards and we leave the proofs in the Appendix \ref{proof Lemma 3,1}. 
\begin{Lem}\label{Lem:Lipschitz norms}
Let $f$ be a function on $\mathbb R$,
$\alpha,\beta\in (0,1]$ and $p, q \ge
0$.
\begin{enumerate}[1)]
\item If $p\le q$, then
$\|f\|_{\alpha,p}<+\infty$ implies
$\|f\|_{\alpha,q}<+\infty$.
\item If $\alpha\le\beta$, then
$\|f\|_{\beta,p}<+\infty$ implies
$\|f\|_{\alpha,p+\beta-\alpha}<+\infty$.
\item If $P$ is a polynomial of degree $d$, then $\|f\|_{\alpha,p}<+\infty$
implies $\|Pf\|_{\alpha,p+d}<+\infty$.
\item Assume that $F$ is a primitive function of $f$,
then $\|f\|_{\alpha,p}<+\infty$ implies
$\|F\|_{1,p+\alpha}<+\infty$. (Hence $\|F\|_{\alpha,p+1}<+\infty$ by 2).)
\end{enumerate}
\end{Lem}

Inspired by \cite{Ba1986}, we introduce the following function space.
\begin{Def}\label{Def:espace des fonctions}
Let $N\ge 0$ be an integer, and $\alpha\in(0,1]$, $p\ge 0$ be two
real numbers. Denote by
$\mathcal H_{\alpha,p}^N$ the vector space of all Borel
functions $h$ on $\mathbb R$ verifying the following
conditions:
\begin{enumerate}[a)]
\item $h$ has $N^{\mathrm{th}}$ order derivative
which is locally of finite variation and which has
finitely many jumps,
\item the continuous part of
$h^{(N)}$ satisfies $\|h_c^{(N)}\|_{\alpha,p}<+\infty$.
%,$h_c^{(N)}$ is absolutely continuous and the
%Radon-Nikodym derivative of $h_c^{(N)}$ w.r.t. the Lebegue measure is $O(|x|^{\alpha+p-1})$ when
%$x\rightarrow\infty$.
\end{enumerate}
\end{Def}

Condition a) implies that the pure jump part of
$h^{(N)}$ is bounded. Condition b) implies that $h^{(N)}_c$ has at most polynomial
increasing speed at infinity, therefore also is $h$.
These conditions allow us to include some irregular  functions such as
indicator functions. Let $k$ be a real
number and $I_k(x)=\indic_{\{x\leq k\}}$. Then $\|I_k\|_{\alpha,p}$ is clearly not finite. However, $\|I_{k,c}\|_{\alpha,p}=0$, which means that
for any $\alpha\in (0,1]$ and any $p\geq 0$, $I_k(x)\in\mathcal
H_{\alpha,p}^{0}$.
%Similarly, $\forall\alpha\in(0,1]$,
%$|x|^{N}\in\mathcal H_{\alpha,0}^N$.
Note that any
function $h$ in $\mathcal H^0_{\alpha,p} $ can be
decomposed as $h=h_c+h_d$,
where $h_c$  satisfies %is absolutely continuous,
$\|h_c\|_{\alpha,p}<+\infty$, 
%and $h_c'(x)=O(|x|^{\alpha+p-1})$, 
the discontinuous part $h_d$ is a linear
combination of indicator functions of the form
$\indic_{\{x\leq k\}}$ plus a constant (so that $h_c(0)=0$).

\begin{Pro}\label{Pro:properties of H}
Let $N\geq 0$ be an integer, $\alpha,\beta\in (0,1]$ and $p, q \ge
0$ be real numbers.
Then the following assertions hold:
\begin{enumerate}[1)]
\item when $N\ge 1$, $h\in\mathcal H_{\alpha,p}^N$ if and only if $h'\in\mathcal
H^{N-1}_{\alpha,p}$;
\item if $p\le q$, then $\mathcal
H_{\alpha,p}^N \subset\mathcal H_{\alpha,q}^N$;
if $\alpha\le\beta$, then $\mathcal H_{\beta,p}^N\subset
\mathcal H_{\alpha,p+\beta-\alpha}^N$;
\item when $N\ge 1$, $\mathcal H_{\alpha,p}^N
\subset\mathcal H_{1,\alpha+p}^{N-1}\subset\mathcal H_{\alpha,p+1}^{N-1}$;
\item if $h\in\mathcal H_{\alpha,p}^N$ and if $P$ is
a polynomial of degree $d$, then $Ph\in\mathcal
H_{\alpha,p+d}^N$.
\end{enumerate}
\end{Pro}
\begin{proof}
1) results from the definition. 2), 3) and 4) are consequences
of Lemma \ref{Lem:Lipschitz norms}.
%4) By 2), it suffices to prove that $x^dh(x)\in\mathcal
%H_{\alpha,p+d}^N $ for any $d\ge 1$. When $N=0$,
%the assertion is true by discussing respectively the two cases when $h$ is continuous and when
%$h$ is an indicator function.
%the function $h-h_c$ is a linear combination of
%indicator functions. Therefore we only need to prove
%the assertion for the case where $h$ is continuous and
%$h$ is an indicator function of the form
%$\indic_{(-\infty,k]}$ respectively. The continuous
%case results from 3) of Lemma \ref{Lem:Lipschitz
%norms}. If $h=\indic_{(-\infty,k]}$, the continuous
%part of $x^dh(x)$ is a primitive function of
%$dx^{d-1}h(x)$. As $\|x^{d-1}h(x)\|_{0,d-1}<+\infty$,
%we obtain $\|(x^dh(x))_c\|_{1,d-1}<+\infty$, and
%therefore $x^dh(x)\in\mathcal
%H^0_{1,d-1}\subset\mathcal H^0_{\alpha,p+d}$. Thus we
%have proved the case where $N=0$.
%The general case
%follows by using the equality
%$(x^dh(x))'=dx^{d-1}h(x)+x^dh'(x)$.
\end{proof}

The following result on the operator $h\rightarrow f_h$
is fundamental. It shows that compared to $h$, the solution of
Stein's equation $f_h$ has one more order in regularity
and its derivative has the same order in increasing speed at infinity. The
proof of this proposition, which is rather technical, is postponed to Appendix \ref{proof Lemma 3.4}.
\begin{Pro}\label{Lem:fundamental lemma}
Assume that $h\in\mathcal H_{\alpha,p}^N$. Then
$f_h\in\mathcal H_{\alpha,p}^{N+1}$.
\end{Pro}

We now restate Theorem \ref{Thm:main theorem}
in the function space context.

\begin{Thm}\label{Thm:main theroem, function space version}
Let $N\ge 0$ be an integer, $\alpha\in(0,1]$ and $p\ge 0$. Assume that  $h\in\mathcal H_{\alpha,p}^N$.  Let
$X_1,\cdots,X_n$ be zero-mean random variables which
have $(N+\max(\alpha+p,2))^{\mathrm{th}}$ order moment. Then
all terms in
\eqref{Equ:recursive CN} and \eqref{Equ:recursive eN}
are well defined, and the equality $\esp[h(W)]=C_N(h)+e_N(h)$ holds.
\end{Thm}
\begin{proof}
When $N=0$, $h\in\mathcal H_{\alpha,p}^0$ and then
 $h(x)=O(|x|^{\alpha+p})$. Hence $\esp[h(W)]$ and
$\Phi_{\sigma_W}(h)$ are well defined.
Assume that we have proved the theorem for
$0,\cdots,N-1$. Let $h\in\mathcal H_{\alpha,p}^N$.
Then by Proposition \ref{Pro:properties
of H},
$h(x)\in\mathcal H_{\alpha,p}^{N}\subset \mathcal H_{\alpha,p+1}^{N-1}\cdots\subset\mathcal
H_{\alpha,p+N}^0$, so $h(x)=O(|x|^{\alpha+p+N})$.
By Proposition \ref{Lem:fundamental lemma}, $f_h\in\mathcal H_{\alpha,p}^{N+1}$ and by Proposition
\ref{Pro:properties of H} 1), for any $|\mathbf J|=1,\cdots,N$, $f_h^{(|\mathbf J|+1)}\in\mathcal
H_{\alpha,p}^{N-|\mathbf J|}$. So the induction hypothesis implies
that $C_{N-|\mathbf J|}(f_h^{(|\mathbf J|+1)})$ and
$e_{N-|\mathbf J|}(f_h^{(|\mathbf J|+1)})$ exist.
Furthermore,  for the terms
$\varepsilon_{N-k}$ and $\delta_N$ in
\eqref{Equ:recursive eN}, since $f_h^{(k+1)}(x)=
O(|x|^{\alpha+p+N-k})$ for any $k=0,\cdots,N$, they are well defined. Finally,
combined with the equality 
$$\esp[(X_i^*)^k]=\frac{\esp[X_i^{k+2}]}{\sigma_i^2(k+1)},$$
all moments figuring in \eqref{Equ:recursive CN} and \eqref{Equ:recursive eN}
exist. Thus all terms  are well defined, and the
formal proof in the previous section shows that
$\esp[h(W)]=C_N(h)+e_N(h)$.
\end{proof}

%==================================
\section{Error estimations}
%==================================

%==========================
\subsection{Concentration inequalities}
%==========================
\hskip\parindent We shall prove some concentration inequalities 
similar to several results in \cite{CS2001,CS2005}, which
give upper bounds for probabilities of the form $\mathbb P(a\leqslant
W\leqslant b)$ with $a$ and $b$ being two real numbers.
We shall take into consideration the parameter $\alpha$ and give some variants where appear
certain lower order moments if $\alpha<1$. When $\alpha=1$, we recover some estimations in \cite{EJ2008}.
These concentration inequalities will be
useful to estimate the approximation error terms
and
the proof is based on the zero bias transformation.

\begin{Lem}\label{Lem:difference X et X*}
Let  $\alpha\in(0,1]$ be a real number and $X$
be a r.v. with mean zero, finite variance
$\sigma^2>0$ and up to $(\alpha+2)^{\mathrm{th}}$ order moments.
Let $X^*$ have the zero biased
distribution of $X$ and be independent of $X$. Then,
for any $\varepsilon>0$,
\[\proba(|X-X^*|>\varepsilon)\le\frac{1}{2\varepsilon^\alpha(\alpha+1)\sigma^2}
\esp[|X^s|^{\alpha+2}],\] where
$X^s=X-\widetilde X$ and $\widetilde X$ is an
independent copy of $X$.
\end{Lem}
\begin{proof}
Similar to the Markov inequality, the following
inequality holds:
\[\proba(|X-X^*|>\varepsilon)\le\frac{1}{\varepsilon^\alpha}
\esp[|X-X^*|^\alpha].\] Moreover, since $X$ and $X^*$ are
independent, the definition of the zero bias transformation (see \cite[Pro2.3]{EJ2008})  implies that
\[\esp[|X-X^*|^\alpha]=\frac{1}{2(\alpha+1)\sigma^2}\esp[|X^s|
^{\alpha+2}].\]
\end{proof}

\begin{Pro}\label{Pro:inequality concentration}Let $X_i \,(i=1,\cdots,n)$ be independent random variables with mean zero and variance $\sigma_i^2>0$. Let $W=X_1+\cdots+X_n$ and denote its variance by $\sigma_W^2$.
For $a,b\in\R$ such that $a\le b$ and any real number $\alpha\in(0,1]$, we have
\begin{equation}\label{Equ:concentration inequality}
\proba(a\le W\le b)\le2\Big(\frac{b-a}{2\sigma_W}\Big)^\alpha+\frac{2}{\alpha+1}
\sum_{i=1}^n\esp\Big[\big|\frac{X_i^s}{\sigma_W}\big|^{\alpha+2}\Big]+\frac{1}{2\sigma_W^2}
\Big(\sum_{i=1}^n\sigma_i^4\Big)^{\frac 12}.\end{equation}
\end{Pro}
\begin{proof}
Let $I_{[a,b]}(x)=1$ if $x\in[a, b]$ and $I_{[a,b]}(x)=0$ otherwise. Its primitive function $f(x):=\int_{(a+b)/2}^{x}I_{[a,b]}(t)dt$ satisfies $|f(x)|\leq (b-a)/2$. Then $$\esp(I_{[a,b]}(W^*))=\frac{1}{\sigma_W^2}\esp(Wf(W))\leq \min\left(\frac{b-a}{2\sigma_W},1\right).$$ Note that for any $u\geq 0$ and any $\alpha\in(0,1]$, $\min(u,1)\leq u^{\alpha}$. Then for any $\varepsilon >0$, 
$$\proba(a-\varepsilon\leq W^*\leq b+\varepsilon)\leq \left(\frac{b-a+2\varepsilon}{2\sigma_W}\right)^\alpha\leq \left(\frac{b-a}{2\sigma_W}\right)^\alpha+\left(\frac{\varepsilon}{\sigma_W}\right)^\alpha$$
where the last inequality is because for any $u$ and $v$ positive, 
one always has $(u+v)^\alpha\leq u^\alpha+v^\alpha$. On the other hand, by using a conditional expectation technique, 
\[\begin{split}\proba(a-\varepsilon\le W^*\le
b+\varepsilon)\ge\proba(a\leq W\leq b,| X_I-X_I^*|\leq\varepsilon)\\
\geq\proba(a\le W\le
b)\proba(|X_I^*-X_I|\le\varepsilon)-\frac14\Big(\sum_{i=1}^n\frac{\sigma_i^4}{
\sigma_W^4}\Big)^{\frac 12}.\end{split}\] 
We recall that $W^*=W^{(I)}+X_I^*$ where $I$ is a random variable taking values in $\{1,\cdots,n\}$ with $\proba(I=i)=\sigma_i^2/\sigma_W^2$,  $W^{(i)}:=W-X_i$, $X_i^*$ has the zero biased distribution of $X_i$ and is independent of $W^{(i)}$. In this proof exceptionally, we assume that $X_i^*$ is also independent of $X_i$. By Lemma \ref{Lem:difference X et X*},  
$$\proba(|X_I^*-X_I|\le\varepsilon)=1-\sum_{i=1}^n
\frac{\sigma_i^2}{\sigma_W^2}\proba({|X_i^*-X_i|>\varepsilon})
\ge1-\frac{1}{2\sigma_W^2(\alpha+1)\varepsilon^\alpha}
\sum_{i=1}^n\esp[|X_i^s|^{\alpha+2}].$$ Finally, the inequality \eqref{Equ:concentration inequality} follows
by taking
\[\varepsilon=\Big(\frac{1}{\sigma_W^2(\alpha+1)}\sum_{i=1}^n
\esp[|X_i^s|^{\alpha+2}]\Big)^{\frac{1}{\alpha}}.\] 

\def\skip{
We shall follow the proof of \cite[Pro3.3]{EJ2008}. Let $\varepsilon>0$ be an arbitrary real number, then
$$\proba(a-\varepsilon\le W^*\le
b+\varepsilon)\le\frac{1}{\sigma_W}\Big(\varepsilon+\frac{b-a}{2}\Big).$$
On the other hand,$$\proba(a-\varepsilon\le W^*\le
b+\varepsilon)\ge\proba(a\le W\le
b)\proba(|X_I^*-X_I|\le\varepsilon)-\frac14\Big(\sum_{i=1}^n\frac{\sigma_i^4}{
\sigma_W^4}\Big)^{\frac 12},$$ where
$$\proba(|X_I^*-X_I|\le\varepsilon)=1-\sum_{i=1}^n
\frac{\sigma_i^2}{\sigma_W^2}\proba({|X_i^*-X_i|>\varepsilon})
\ge1-\frac{1}{2\sigma_W^2(\alpha+1)\varepsilon^\alpha}
\sum_{i=1}^n\esp[|X_i^s|^{\alpha+2}].$$
The inequality \eqref{Equ:concentration inequality} follows
by taking
\[\varepsilon=\Big(\frac{1}{\sigma_W^2(\alpha+1)}\sum_{i=1}^n
\esp[|X_i^s|^{\alpha+2}]\Big)^{\frac{1}{\alpha}}.\]}
\end{proof}

\begin{Cor}\label{Cor:concentration inequatliyt for Wi}
Let $a,b\in\R$ such that $a\le b$ and $i\in\{1,\cdots,n\}$, then
\[\begin{split}\proba(a\le W^{(i)}\le b)\le 4\Big(\frac{b-a}{2\sigma_W}\Big)^\alpha+\frac{4}{\alpha+1}\sum_{j=1}^n\esp\Big[\big|\frac{X_j^s}{\sigma_W}\big|^{\alpha+2}\Big]
+\frac{1}{\sigma_W^2}\big(\sum_{j=1}^n{\sigma_j}^4\big)^{\frac
12}+4\Big(\frac{2\sigma_i}{\sigma_W}\Big)^\alpha
\end{split}\]where $W^{(i)}=W-X_i$.
\end{Cor}
\begin{proof}
Let $\varepsilon>0$ be a real number, then
\[\mathbb P(a\le W^{(i)}\le b,\;|X_i|\le\varepsilon)\le
\mathbb P(a-\varepsilon\le W\le b+\varepsilon).\] Note
that $W^{(i)}$ and $X_i$ are independent and
\[\mathbb P(|X_i|\le\varepsilon)=1-\mathbb P(|X_i|>\varepsilon)
\ge 1-\frac{E[|X_i|]}{\varepsilon}.\] By Proposition \ref{Pro:inequality concentration} and taking
$\varepsilon=2\mathbb E[|X_i|]$, we obtain the
inequality.
\end{proof}

%=======================================
\subsection{Estimations of error terms}
%=======================================

In this section, we shall estimate the
error term $e_N(h)$ in Theorem \ref{Thm:main theroem,
function space version}. The recursive formulas \eqref{Equ:recursive eN}
and \eqref{equ:reverse taylor remainder} permit us to reduce the
problem to the estimation of classical Taylor expansion
errors. %We shall keep the notation of Theorem
%\ref{Thm:main theorem}.

For any positive random variable $Y$ and any real number $\beta\ge 0$,
we introduce the notation $m_{Y}^{(\beta)}:=
\mathbb E[Y^\beta]/\Gamma(\beta+1)$, where $\Gamma$ is
the Gamma function $\Gamma(x)=\int_0^\infty t^{x-1}e^{-t}dt$. This
notation generalizes the one introduced in \S
\ref{Sec:introduction} since when $\beta\in\mathbb N$,
$\Gamma(\beta+1)=\beta!$. 
\begin{Pro}\label{Pro:remaining term of Taylor}
Let $N\ge 0$ be an integer, $\alpha\in(0,1]$ and $p\ge 0$ be two
real numbers. Let $X$ be a random
variable which has  up to $(N+\alpha+p)^{\mathrm{th}}$ moments and
satisfies the following concentration inequality
\[\proba(a\le X\le b)\le c(b-a)^{\alpha}+r,\qquad
\forall a,b\in\mathbb R, \,a\leq b,\] where $c$ and $r$ are two constants. Let $Y$ be a random variable
which is independent of $X$ and has up to
$(N+\alpha+p)^{\mathrm{th}}$ moments. Then, for
any function $g\in\mathcal H_{\alpha,p}^N$ and any
$k=0,\cdots,N$,
\begin{equation}\label{Equ:estime of delta}\begin{split}
\big|\delta_{N-k}(g^{(k)},X,Y)\big|
&\le V(g^{(N)}_d)\Big(cm_{|Y|}^{(N-k+\alpha)}+r
m_{|Y|}^{(N-k)}\Big)\\&
\;\;+\|g^{(N)}_c\|_{\alpha,p}
\Big(u_{\alpha,p,X}m_{|Y|}^{(N-k+\alpha)}+v_{\alpha,p}
m_{|Y|}^{(N-k+\alpha+p)}\Big),
\end{split}\end{equation}
where $V(g^{(N)}_d)$ denotes the total variation of
$g^{(N)}_d$, the coefficients $u_{\alpha,p,X}$ and $v_{\alpha,p}$ are defined as
$u_{\alpha,p,X}=\big(1+(1+2^p)\esp[|X|^p]\big)\Gamma(\alpha+1)$ and
$v_{\alpha,p}=2^p\Gamma(\alpha+p+1)$.
\end{Pro}

\begin{proof}
We have by \eqref{Equ:remaining term of Taylor expansion} that when $k<N$,
\[\delta_{N-k}(g^{(k)},X,Y)=\frac{1}{(N-k-1)!}
\int_0^1(1-t)^{N-k-1}\esp\big[
\big(g^{(N)}(X+tY)-g^{(N)}(X)\big)Y^{N-k}
\big]\,dt.\] Since $g\in\mathcal
H_{\alpha,p}^{N}$, the function $g^{(N)}_d$ can be written as
$$g^{(N)}_d(x)=g_d^{(N)}(0)+\sum_{1\leq j\leq M}\varepsilon_j\indic_{x\leq K_j}-\sum_{\begin{subarray}{c}1\leq j\leq M\\K_j\geq 0\end{subarray}}\varepsilon_j.$$
Therefore,
$g^{(N)}_d(X+tY)-g^{(N)}_d(X)=\sum_{j=1}^M\varepsilon_j
\indic_{K_j-tY_+<X\le K_j-tY_-}$, where
$Y_+=\max(Y,0)$ and $Y_-=\min(Y,0)$. Thus the
concentration inequality hypothesis implies that
\[\esp\big[|g^{(N)}_d(X+tY)-g^{(N)}_d(X)|\;\big|Y\big]\le
\sum_{j=1}^M|\varepsilon_j|\big(c t^\alpha|Y|^\alpha+r\big).\]
Moreover, one has $$\int_0^1\frac{(1-t)^{N-k-1}}{(N-k-1)!}\esp\Big[\sum_{j=1}^M|\varepsilon_j|\big(c t^\alpha|Y|^\alpha+r\big)|Y|^{N-k}\Big]dt
=\sum_{j=1}^M|\varepsilon_j|\big(cm_{|Y|}^{(N-k+\alpha)}+r
m_{|Y|}^{(N-k)}\big)$$ by using the
following equality concerning Beta function
\[B(x,y):=\int_0^1t^{x-1}(1-t)^{y-1}\,dt=\frac{\Gamma(x)\Gamma(y)}{\Gamma(x+y)}, \quad x, y>0.\]
On the other hand, by definition of the norm
$\|\cdot\|_{\alpha,p}$, we have
\[\begin{split}\big|g^{(N)}_c(X+tY)-g^{(N)}_c(X)\big|
&\le\|g_c^{(N)}\|_{\alpha,p}|tY|^\alpha(1+|X+tY|^p+|X|^p)\\
&\le\|g^{(N)}_c\|_{\alpha,p}
|tY|^\alpha\big(1+(2^p+1)|X|^p+2^p|tY|^p\big),\end{split}\]
%\[\begin{split}\esp\Big[\big|g^{(N)}_c(X+tY)-g^{(N)}_c(X)\big|\;\Big|\;Y\Big]
%&\le\|g_c^{(N)}\|_{\alpha,p}|tY|^\alpha(1+|X+tY|^p+|X|^p)\\
%&\le\|g^{(N)}_c\|_{\alpha,p}
%|tY|^\alpha\big(1+(2^p+1)\esp[|X|^p]+2^p|tY|^p\big),\end{split}\]
where the last inequality results from $(a+b)^p\le
2^p(a^p+b^p)$. Note that
$$\int_0^1\frac{(1-t)^{N-k-1}}{(N-k-1)!}\esp\Big[|tY|^\alpha\big(1+(2^p+1)|X|^p+2^p|tY|^p\big)|Y|^{N-k}\Big]dt
=u_{\alpha,p,X}m_{|Y|}^{(N-k+\alpha)}+v_{\alpha,p}
m_{|Y|}^{(N-k+\alpha+p)}.$$
Thus we obtain the estimation \eqref{Equ:estime of delta}. \\Finally, it remains to check the case when $k=N$.
Consider the continuous and discontinuous parts of
$\delta_{0}(g^{(N)},X,Y)=\esp\big[g^{(N)}(X+Y)-g^{(N)}(X)\big]$ respectively.
By using similar method as above, we obtain that
$\esp\big[|g^{(N)}_d(X+Y)-g^{(N)}_d(X)|\big]\le
V(g^{(N)}_d)\big(c m_{|Y|}^{(1)}+r\big)$ and $\esp\big[|g^{(N)}_c(X+Y)-g^{(N)}_c(X)|\big]\le
\|g^{(N)}_c\|_{\alpha,p}\left( \esp[|Y|^{\alpha}]\big(1+(2^p+1)\esp[|X|^p]\big)+2^p\esp[|Y|^{\alpha+p}]\right)$,
which implies \eqref{Equ:estime of delta}.
\end{proof}

By Proposition \ref{Lem:reverse taylor formula} and Proposition \ref{Pro:remaining term of Taylor},
we obtain the error estimation for the reverse Taylor expansion.
\begin{Cor}\label{Cor:remaining term of Taylor inverse}
With the previous notation, we have
\begin{equation}\label{equ: epsilon N}\begin{split}|\varepsilon_N(g,X,Y)|\le\sum_{d\ge 0}
\hspace{-0mm}\sum_{\mathbf{J}\in\mathbb N^d_*,\,|\mathbf{J}|\le N}
%\hspace{-4mm}
&m_{|Y|}^{(\mathbf{J})}\Big[V(g_d^{(N)})\big(cm_{|Y|}^{(N-|\mathbf{J}|+\alpha)}
+rm_{|Y|}^{(N-|\mathbf{J}|)}\big)\\
&+\|g_c^{(N)}\|_{\alpha,p}\big(u_{\alpha,p,X}m_{|Y|}^{(N-|\mathbf{J}|+\alpha)}
+v_{\alpha,p}m_{|Y|}^{(N-|\mathbf{J}|+\alpha+p)}\big)\Big],
\end{split}\end{equation}
\end{Cor}

Combining the concentration inequality (Corollary \ref{Cor:concentration inequatliyt for Wi}) 
and the above estimations (Proposition \ref{Pro:remaining term of Taylor} and
Corollary \ref{Cor:remaining term of Taylor inverse}), we obtain upper bounds for the Taylor and reverse Taylor remainders  $\delta_{N-k}(f_h^{(k)},W^{(i)}, X_i)$ and $\varepsilon_{N-k}(f_h^{(k)},W^{(i)}, X_i)$, where the summand variables $X_1,\cdots, X_n$ are independent. 
This allows us, together with the  recursive formula \eqref{Equ:recursive eN}, 
to obtain an upper bound for the asymptotic expansion remainder $e_N(h)$.

In particular, we give in the following the order estimation of $e_N(h)$ when  
$X_1,\cdots,X_n$ are in addition  i.i.d. random variables.

\begin{Pro}\label{Cor:estimation de delta et eps}
Suppose that $X_1, \cdots,X_n$
are i.i.d. random variables with mean zero and up to $(N+2+\alpha+p)^{\text{th}}$ order moments, normalized such that each $X_i$ has the same distribution as $Z/\sqrt n$ where $Z$ is a fixed random variable with mean zero and finite non-zero variance.
Then for any function $g\in\mathcal
H_{\alpha,p}^N$ and any $k=0,\cdots,N$, we have
\begin{eqnarray}\label{equ:error order}
\delta_{N-k}(g^{(k)},W^{(i)},X_i^*)=O \left(\Big(\frac{1}{\sqrt{n}}\Big)^{N-k+\alpha+p}\right),\\
\label{equ:error order2}
\varepsilon_{N-k}(g^{(k)},W^{(i)},X_i)=O\left(
\Big(\frac{1}{\sqrt{n}}\Big)^{N-k+\alpha+p}\right),
\end{eqnarray}
where $W^{(i)}=W-X_i$ and $X_i^*$ is independent of $W^{(i)}$. The
implied constants depend on $\|g^{(N)}_c\|_{\alpha,p}$,
$V(g^{(N)}_d)$ and up to
$(N-k+2+\alpha+p)^{\text{th}}$ order moments of
$Z$.
\end{Pro}
\begin{proof}
By Corollary
\ref{Cor:concentration inequatliyt for Wi}, we have for any $a\leq b$ and any $\alpha\in
(0,1]$ that
\[\proba(a\le W^{(i)}\le b)\le c(b-a)^\alpha+r(n)\]
where the coefficients are given by
\[c=\frac{2^{2-\alpha}}{\sigma^\alpha},\qquad
r(n)=\frac{4}{\sigma^{2+\alpha}(\alpha+1)}\frac{\esp[|Z^s|^{\alpha+2}]}{\sqrt{ n}^\alpha}
+\frac{1}{\sqrt n}+\frac{8}{\sqrt{n}^\alpha}.\]
By Proposition \ref{Pro:remaining term of Taylor}, we obtain an upper bound of
$\delta_{N-k}(g^{(k)},W^{(i)},X_i^*)$ which is determined by a linear combination of terms (with coefficient not depending on $n$): 
\begin{equation} \label{equ: coefficients}m_{|X_i^*|}^{(N-k+\alpha)}, \,\,r(n)m_{|X_i^*|}^{(N-k)},\,\,  
\esp[|W^{(i)}|^p]m_{|X_i^*|}^{(N-k+\alpha)}\,\, \text{ and }\,\,
m_{|X_i^*|}^{(N-k+\alpha+p)}.\end{equation} Note that $r(n)=O((1/\sqrt{n})^\alpha)$.  For any $k=0,\cdots,N$, $\esp[|X_i^*|^k]$ equals
${\esp[|X_i|^{k+2}]}/{(\sigma_i^2(k+1))}$ and is of order $(1/\sqrt{n})^k$.
So the first three terms in \eqref{equ: coefficients} are of
order $(1/\sqrt{n})^{N-k+\alpha}$ and the last term is of order
$(1/\sqrt{n})^{N-k+\alpha+p}$, which implies the first assertion.
The second assertion then follows by Corollary \ref{Cor:remaining term of Taylor inverse}.

\def\skip{
Note that $r\ll (1/\sqrt{n})^\alpha$
(when $\alpha>0$, one even has $r\ll 1/\sqrt{n}$).
Furthermore, for any real number $\beta\ge 0$, one has
\begin{equation}\label{Equ:moment de xistar}m_{|X_i|}^{(\beta)}=\Big(\frac{1}{\sqrt{n}}\Big)^\beta
m_{|Z|}^{(\beta)}\quad\;\text{and}\quad\;
m_{|X_i^*|}^{(\beta)}=\sigma_i^{-2}(\beta+2)m_{|X_i|}^{(\beta+2)}
=\Big(\frac{1}{\sqrt{n}}\Big)^\beta
\sigma^{-2}(\beta+2)m_{|Z|}^{(\beta+2)}.\end{equation}}
\end{proof}

\begin{Rem}
According to \eqref{Equ:estime of delta} and \eqref{equ: epsilon N}, the implicit constants in \eqref{equ:error order} and \eqref{equ:error order2}
can be explicitly calculated.
\end{Rem}
\begin{Pro}\label{Pro:error estimation rigorous one}
Let $N\ge 0$ be an integer, $\alpha\in(0,1]$ and $p\ge 0$ be two
real numbers. Let $h$ be a
function in $\mathcal H_{\alpha,p}^N$, and
$X_1,\cdots,X_n$ be as in Proposition \ref{Cor:estimation de delta et eps}.
Then the error term $e_N(h)$ defined in \eqref{Equ:recursive eN} satisfies $$e_N(h)=O
\left(\Big(\frac{1}{\sqrt{n}}\Big)^{N+\alpha+p}\right),$$
where the implied constant depends on up to
$(N+2+\alpha+p)^{\mathrm{th}}$ order moment of $Z$.
\end{Pro}
\begin{proof}
We prove the theorem by induction on $N$. When $N=0$,
\[e_0(h)=\sum_{i=1}^n\sigma_i^2\Big(\delta_0(f_h',W^{(i)},X_i^*)+\varepsilon_0(f_h',W^{(i)},X_i)\Big).\]
Since $h\in\mathcal H_{\alpha,p}^0$, $f_h\in\mathcal H_{\alpha,p}^1$. Then by
Proposition \ref{Cor:estimation de delta et eps},
$e_0(h)=O(({1}/{\sqrt{n}})^{\alpha+p})$.
\def\skip{
\begin{equation}\label{Equ:recursive eN}\begin{split}
e_N(h)&=\sum_{i=1}^n\sigma_i^2\bigg[ \sum_{d\ge
1}(-1)^{d-1}\hspace{-5mm}\sum_{\mathbf{J}\in\mathbb
N_*^d,\,|\mathbf{J}|\le
N}\hspace{-4mm}m_{X_i}^{(\mathbf
J^\circ)}\big(m_{X_i^*}^{(\mathbf
J^\dagger)}-m_{X_i}^{(\mathbf{J}^\dagger)}\big)e_{N-|\mathbf{J}|}(f_h^{(|\mathbf{J}|+1)})\\&
\qquad
+\sum_{k=0}^N\varepsilon_{N-k}(f_h^{(k+1)},W^{(i)},X_i)m_{X_i^*}^{(k)}
+\delta_N(f_h',W^{(i)},X_i^*)\bigg],
\end{split}\end{equation}
}

\noindent Assume that we have already proved the theorem for
$0,\cdots, N-1$. Consider $h\in\mathcal H^N_{\alpha,p}$ and $e_N(h)$ defined as in \eqref{Equ:recursive eN}.
For any $\mathbf J$ such that $1\leq |\mathbf{J}|\leq N$, $e_{N-|\mathbf{J}|}=O(\displaystyle
({1}/{\sqrt{n}})^{N-|\mathbf{J}|+\alpha+p})$. In addition, since
$|\mathbf J^\circ|+|\mathbf J^\dagger|=|\mathbf J|$, we have that $m_{X_i}^{(\mathbf
J^\circ)}\big(m_{X_i^*}^{(\mathbf
J^\dagger)}-m_{X_i}^{(\mathbf{J}^\dagger)}\big)$ is of order $({1}/{\sqrt{n}})^{|\mathbf J|}$.
On the other hand, $f_h'\in\mathcal H^{N}_{\alpha,p}$, so $\delta_N(f_h',W^{(i)},X_i^*)
=O(\displaystyle({1}/{\sqrt{n}})^{N+\alpha+p})$.
Moreover, for any $k=0,\cdots,N$, $f_h^{(k+1)}\in\mathcal H^{N-k}_{\alpha,p}$. So
$\varepsilon_{N-k}(f_h^{(k+1)},W^{(i)},X_i^*)
=O(\displaystyle({1}/{\sqrt{n}})^{N-k+\alpha+p})$. Finally we have $m_{X_i^*}^{(k)}=O
(({1}/{\sqrt{n}})^{k})$. Combining all the above estimations, we prove the proposition.
\def\skip{
have already shown that, for any $k\in\{0,\cdots, N\}$,
$f_h^{(k+1)}\in\mathcal H^{N-k}_{\alpha,p}$. Hence the
recurrence hypothesis implies that
$e_{N-k}(f_h^{(k+1)})\ll(\frac{1}{\sqrt{n}})^{N-k+\alpha}$,
$k=1,\cdots,N$, where the constant depends only
on up to $(N-k+\alpha+p+2)^{\mathrm{th}}$ order moments
of $|Z|$. Combining with the estimations in Corollary
\ref{Cor:estimation de delta et eps}, together with the
computation of moments \eqref{Equ:moment de xistar}, we
obtain that
$e_N(h)\ll\displaystyle\Big(\frac{1}{\sqrt{n}}\Big)^{N+\alpha}$,
where the constant depends only on up to
$(N+\alpha+p+2)^{\mathrm{th}}$ order moment.}
\end{proof}

Consider now several examples. Let $I_k(x)=\indic_{\{x\leq k\}}$ be the indicator function. As mentioned before,
$I_k\in\mathcal H^0_{\alpha,0}$. By Proposition
\ref{Pro:error estimation rigorous one}, we know that if $X_1,\cdots,X_n$ are i.i.d. random variables with up to $(2+\alpha)^{\text{th}}$ order moment, then $e_0(h)=O((1/\sqrt{n})^\alpha)$, where
 the coefficient depends on up to $(2+\alpha)^{\text{th}}$ moment of the summand variables. 
 This is similar to a result (Theorem 6) in \cite[\S V.3]{Pe1975}. When $\alpha=1$, it corresponds to the order estimation in the classical Berry-Esseen inequality.

Let $h(x)=(x-k)_+$ be the call function discussed in \cite{EJ2008}. As a primitive function of the indicator function,
we know that $h\in\mathcal H^1_{\alpha,0}$. So the call function admits a first order expansion given
by \eqref{Equ:recursive CN} as:
$$C_1(h)=\Phi_{\sigma_W}(h)+\sum_{i=1}^n\sigma_i^2\esp[X_i^*]\Phi_{\sigma_W}(f_h'').$$
Moreover, since $\sigma_W^2\Phi_{\sigma_W}(f_h'')=\Phi_{\sigma_W}(xf_h')
=\frac{1}{\sigma_W^2}\Phi_{\sigma_W}\big((\frac{x^2}{3\sigma_W^2}-1)xh(x)\big)$. We recover the
correction term in \cite{EJ2008}. 

%%%%%%%%%%%%%%%%%%%%%%%%%%%
\appendix
%%%%%%%%%%%%%%%%%%%%%%%%%%%

%%%%%%%%%%%%%%%%%%%%%%%%%%%%
\section{Proof of Lemma \ref{Lem:Lipschitz norms}}\label{proof Lemma 3,1}
%%%%%%%%%%%%%%%%%%%%%%%%%%%%

\begin{proof} For the first two assertions, it
suffices to prove respectively the boundness
of the following two functions
\[\frac{1+|x|^p+|y|^p}{1+|x|^q+|y|^q},\qquad
|x-y|^{\beta-\alpha}\frac{1+|x|^p+|y|^p}
{1+|x|^{p+\beta-\alpha}
+|y|^{p+\beta-\alpha}}.\] These functions are
both continuous on $\mathbb R$, therefore are
bounded on any compact subset of $\mathbb
R^2$. Thus we may assume without loss of generality that $r=\sqrt{x^2+y^2}\ge 1$. In
this case, $\max\{|x|,|y|\}\ge
{r}/{\sqrt{2}}$, so
\begin{gather*}
\frac{1+|x|^p+|y|^p}{1+|x|^q+|y|^q}\le\frac{1+2r^p}{1+(r/\sqrt{2})^q}
\le 3\cdot 2^{q/2},\\
|x-y|^{\beta-\alpha}\frac{1+|x|^p+|y|^p}{1+|x|^{p+\beta-\alpha}
+|y|^{p+\beta-\alpha}}\le(2r)^{\beta-\alpha}\frac{1+2r^p}
{1+(r/\sqrt{2})^{p+\beta-\alpha}}\leq 3\cdot
2^{(p+3\beta-3\alpha)/2}.
\end{gather*}

3) One has
\[\frac{|P(x)f(x)-P(y)f(y)|}{|x-y|^\alpha(1+|x|^{p+d}
+|y|^{p+d})}\le
\frac{(1+|x|^p+|y|^p)P(x)}{1+|x|^{p+d}+|y|^{p+d}}\|f\|_{\alpha,p}
+\frac{|f(y)|\cdot|P(x)-P(y)|}{|x-y|^\alpha(1+|x|^{p+d}+|y|^{p+d})}.\]
By using the argument as in the proof of 1)
and 2), we obtain that the first term in the
right-hand side is bounded. Since $P$ is a
polynomial of degree $d$, there exists a
polynomial $Q(x,y)$ in two variables and of
degree $d-1$, such that
$Q(x,y)=(P(x)-P(y))/(x-y)$. Therefore, the
second term equals
\[\frac{|Q(x,y)|\cdot|x-y|^{1-\alpha}\cdot|f(y)|}{1+|x|^{p+d}+|y|^{p+d}}\]
which is bounded by a similar argument as for
proving 1) and 2).

4) Since $\|f\|_{\alpha,p}<+\infty$,  $|f(t)|\ll 1+|t|^{\alpha+p}$. Therefore, for
$x,y\in\mathbb R$, $x\le y$, one has
\[|F(x)-F(y)|\le\int_x^y|f(t)|\,dt\ll\int_x^y
\big(1+|t|^{\alpha+p}\big)\,dt\le
(1+|x|^{\alpha+p}+|y|^{\alpha+p})|y-x|.\]
Hence
$\displaystyle\frac{|F(x)-F(y)|}{|x-y|(1+|x|^{p+\alpha}+|y|^{p+\alpha})}
$ is bounded.
\end{proof}

%%%%%%%%%%%%%%%%%%%%%%%%%%
\section{Proof of Proposition \ref{Lem:fundamental lemma}}\label{proof Lemma 3.4}
%%%%%%%%%%%%%%%%%%%%%%%%%%

\hskip\parindent We now prove the
Proposition \ref{Lem:fundamental lemma}. Let
$h\in\mathcal H_{\alpha,p}^N$. The function
$f_h$ is one more order differentiable than
$h$ and is hence $N+1$ times
differentiable. Taking $N^\text{th}$ order
derivative on both sides of Stein's equation,
we get
\begin{equation}\label{Equ:derive d'ordre N de Stein equation}
(xf_h(x))^{(N)}-\sigma^2 f^{(N+1)}_h(x)=
h^{(N)}(x).\end{equation}
The function
$(xf_h(x))^{(N)}$ is continuous, so $f_h^{(N+1)}$ is locally of finite variation and has
finitely many jumps as $h^{(N)}(x)$ does. In the following, we shall prove
$\|f_{h,c}^{(N+1)}\|_{\alpha,p}<+\infty$. 

\begin{Def}
Let $A$ be an interval in $\mathbb R$ and $f$
be a Borel function on $A$. For any $\alpha\in(0,1]$ and $p\ge 0$, we define
\begin{equation}\label{Equ: norm f on A}
\|f\|_{\alpha,p}^A:=\sup_{\begin{subarray}{c}x\neq
y\\ x,y\in
A\end{subarray}}\frac{|f(x)-f(y)|}{|x-y|^\alpha(
1+|x|^p+|y|^p)}.\end{equation}
\end{Def}
This definition is analogous to \eqref{equ:f norm alpha p}, restricted to an interval. When $A$ avoids an open neighborhood of $0$, then the finiteness of $\|f\|_{\alpha,p}^A$ is equivalent to that of $\sup_{\begin{subarray}{c}x\neq
y\\ x,y\in A\end{subarray}}\frac{|f(x)-f(y)|}{|x-y|^\alpha(
|x|^p+|y|^p)}$. This property does not hold for the norm $\|.\|_{\alpha,p}$ defined in \eqref{equ:f norm alpha p}. As a consequence, we have the following result.

\begin{Lem}\label{Lem:5.2}
Let $A\subset (-\infty,-1]\cup[1,+\infty)$ be an interval, $\alpha\in (0,1]$ and $p\geqslant 0$. Let
$q$ be a real number such that $0\leqslant
q\leqslant p$. Then for any Borel function
$f$ defined on $A$,
$\|f\|_{\alpha,p}^A<+\infty$ if and only if
$\|f(x)/x^{p-q}\|_{\alpha,q}^A<+\infty$.
\end{Lem}
\begin{proof} If $\|f(x)/x^{p-q}\|_{\alpha,q}^A<+\infty$, then by similar arguments as for proving Lemma \ref{Lem:Lipschitz norms}, we have  $\|f\|_{\alpha,p}^A<+\infty$.
We now consider the converse assertion. 
Firstly, there exists a constant $C>0$ such that $|f(x)|\leqslant C |x|^{\alpha+p}$ for any $x\in A$.
For any $x,y\in A$, $|x|<|y|$, 
\[\frac{|f(x)x^{q-p}-f(y)y^{q-p}|}{|x-y|^\alpha(1+|x|^q+|y|^q)}
\leqslant |f(x)|\frac{|x^{q-p}-y^{q-p}|}{|x-y|^\alpha(1+|x|^q+|y|^q)}
+|y|^{q-p}\frac{|f(x)-f(y)|}{|x-y|^\alpha(1+|x|^q+|y|^q)}.
\]
The second term is finite since  
$$\frac{|f(x)-f(y)|}{|x-y|^\alpha(|y|^{p-q}+|x|^q|y|^{p-q}+|y|^p)}\leq\frac{|f(x)-f(y)|}{|x-y|^\alpha(1+|x|^p+|y|^p)}=\|f\|_{\alpha,p}^A.$$ By the mean value theorem, the first term is bounded by
\[C|x|^{\alpha+p}\frac{|x-y|\cdot|q-p|\cdot|x|^{q-p-1}}{|x-y|^{\alpha}(
1+|x|^q+|y|^q)}\]
and thus by $C|q-p|$ if we assume in addition that $|y|<2|x|$. When $|y|\geqslant 2|x|$, one has \[|f(x)|\frac{|x^{q-p}-y^{q-p}|}{|x-y|^\alpha(1+|x|^q+|y|^q)}
\leqslant C\frac{|x|^{\alpha+p}|x|^{q-p}}{|x|^{\alpha+q}}\leqslant C.\]
\end{proof}

The following lemma allows us to consider the 
estimations on several disjoint intervals
respectively.

\begin{Lem}\label{Lem:two intervals}
If $A=A_1\cup
A_2$ where $A_1$ and $A_2$ are two intervals
such that $A_1\cap A_2\neq\emptyset$, then
\[\sup\{\|f\|_{\alpha,p}^{A_1},\|f\|_{\alpha,p}^{A_2}\}
\le \|f\|_{\alpha,p}^{A}\le
2(\|f\|_{\alpha,p}^{A_1}+\|f\|_{\alpha,p}^{A_2}).\]
\end{Lem}
\begin{proof}
The first inequality is obvious. For the
second inequality, we only need to prove for
any $x\in A_1$ and $y\in A_2$ that
\[\frac{|f(x)-f(y)|}{|x-y|^\alpha( 1+|x|^p+|y|^p)}
\le
2(\|f\|_{\alpha,p}^{A_1}+\|f\|_{\alpha,p}^{A_2}).\]
Without loss of generality, we may suppose
that $A_1\cap A_2$ contains a single point
$z$. Then $|f(x)-f(y)|\leq
|f(x)-f(z)|+|f(y)-f(z)|$. In addition, since
$z$ is between $x$ and $y$, we have
$|x-y|\ge\max(|x-z|,|y-z|)$ and
$|x|^p+|y|^p\ge\frac
12\max(|x|^p+|z|^p,|y|^p+|z|^p)$. So
\[\frac{|f(x)-f(y)|}{|x-y|^\alpha( 1+|x|^p+|y|^p)}
\le 2\left(\frac{|f(x)-f(z)|}{|x-z|^\alpha(
1+|x|^p+|z|^p)}+\frac{|f(z)-f(y)|}{|z-y|^\alpha(
1+|z|^p+|y|^p)}\right),\] which implies the
second inequality.
\end{proof}

\begin{Lem}\label{Pro:bounded on finite interfval}
If $h\in\mathcal H_{\alpha,p}^N$, then
$\|f_{h,c}^{(N+1)}\|_{\alpha,p}^{A}<+\infty$ for any
bounded interval $A$.
\end{Lem}
\begin{proof}
Firstly, for any bounded interval $A$ and any
Borel function $g$,
$\|g\|_{\alpha,p}^A<+\infty$ if and only if
$g$ is $\alpha$-Lipschitz on $A$. We examine
$f_{h,c}^{(N+1)}$ using \eqref{Equ:derive
d'ordre N de Stein equation}. Since
$h\in\mathcal H_{\alpha,p}^N$, $h^{(N)}_c$ is
locally $\alpha$-Lipschitz. The function
$(xf_h(x))^{(N+1)}=xf_h^{(N+1)}(x)+(N+1)f_h^{(N)}(x)$
has finitely many jumps. Hence
$(xf_h(x))^{(N)}$ is a primitive function of
a locally bounded function, thus is locally
$1$-Lipschitz.  So by \eqref{Equ:derive
d'ordre N de Stein equation},
$f_{h,c}^{(N+1)}$ is locally
$\alpha$-Lipschitz, which implies the
lemma.
\end{proof}

Let $A_1=[-1,1]$, $A_2=(-\infty,-1]$
and $A_3=[1,+\infty)$. Lemma \ref{Lem:two
intervals} shows that to prove the finiteness
of $\|f_{h,c}^{(N+1)}\|$, it suffices to
prove respectively the finiteness of
$\|f_{h,c}^{(N+1)}\|_{\alpha,p}^{A_i},\,(i=1,2,3)$.
Lemma \ref{Pro:bounded on finite interfval}
shows that
$\|f_{h,c}^{(N+1)}\|_{\alpha,p}^{[-1,1]}<+\infty$.
So it remains to deal with $f_{h,c}^{(N+1)}$
on the set $A_2\cup A_3=\R\setminus (-1,1)$.
To this end, we introduce a ``modified''
Stein's equation as in
\cite[Appendix]{EJ2008}:
\begin{equation}\label{Equ:modified Stein equatli}
x\widetilde f_h(x)-\sigma^2\widetilde f_h'(x)=h(x),\quad x\in\mathbb R\setminus(-1,1)\end{equation}
whose solution is given by
\begin{equation}\label{Equ:solution modified Stein equation}
\widetilde f_h(x):=\begin{cases}
\frac{1}{\sigma^{2}\phi_\sigma(x)}\int_x^\infty h(t)
\phi_\sigma(t)\,dt,&x\ge 1,\\\frac{1}{\sigma^{2}\phi_\sigma(x)}
\int^x_{-\infty} h(t) \phi_\sigma(t)\,dt&x\le
-1.\end{cases}\end{equation}
Working with \eqref{Equ:solution modified Stein equation},
it will be easier to treat the derivative functions. In fact, in \eqref{Equ:solution of Stein's equation},
the integrand function $h-\Phi_\sigma(h)$ is centralized under the normal expectation. However, it is not
the case when taking derivatives. This is one reason why we introduce \eqref{Equ:modified Stein equatli}.
Note that in general, the right-hand side of \eqref{Equ:solution modified Stein equation} can not be extended as
a continuous function on $\R$, except
in the special case $\Phi_\sigma(h)=0$ where we recover the solution of classical Stein's equation.

To study $\widetilde f_h$, we introduce
the function space $\mathcal E_\sigma$: for any $\sigma>0$,
let $\mathcal E_\sigma$ be the space of all Borel
functions $h$ on $\R\setminus (-1,1)$  such
that $ \int_{|x|\geq 1}
|h(x)P(x)|\phi_{\sigma}(x)dx<\infty$ for any
polynomial $P$. Note that $\mathcal E_\sigma$
is a vector space which contains all Laurent
polynomials (that is, polynomials in $x$ and
$x^{-1}$) and is stable by multiplication by
Laurent polynomials. Furthermore, as shown by the lemma below, it is
invariant by the operator $h\rightarrow \widetilde f_h$.

\begin{Lem}
Let $h\in\mathcal E_\sigma$.  Then the function $\widetilde f_h$ is
well defined and $\widetilde
f_h\in\mathcal E_\sigma$. Furthermore, if $H$ is a primitive function of $h$, then $H\in\mathcal E_\sigma$.
\end{Lem}
\begin{proof}
Let $P$ be an arbitrary polynomial on $\mathbb R$.
Then
\[\int_1^\infty \hspace{-2mm}|P(x)\widetilde f_{h}(x)|\phi_\sigma(x)\,dx
\le
\frac{1}{\sigma^2}\int_1^\infty\hspace{-2mm}dx\,|P(x)|
\int_x^\infty
\hspace{-2mm}|h(t)|\phi_\sigma(t)\,dt=\frac{1}{\sigma^2}
\int_1^\infty\hspace{-2mm}dt\,|h(t)|\phi_\sigma(t)
\int_1^t|P(x)|\,dx.\] There exists a
polynomial $Q$ such that
$\int_1^t|P(x)|\,dx\le Q(t)$ for any $t\ge 0$.
Therefore, the fact that $h\in\mathcal E_\sigma$
implies that $\int_1^\infty|P(x)\widetilde
f_h(x)|\phi_\sigma(x)\,dx<+\infty$. The
finiteness of the integral on $(-\infty,1]$ is similar. 
The second assertion can be proved by integration by part.
\end{proof}

\begin{Rem}
Note that $\widetilde f_h$ is the only solution of \eqref{Equ:modified Stein equatli} in $\mathcal E_\sigma$, provided that $h\in\mathcal E_\sigma$.
\end{Rem}

More generally, for the derivatives of
$\widetilde f_h$, we consider, for any
integer $N\geq 1$, the set $\mathcal
E_\sigma^N$ which contains all functions $h$
such that $h$ is $N$ times
differentiable on $ \R\setminus (-1,1)$ and
that $h^{(N)}\in\mathcal E_\sigma.$ It is not difficult to
observe that  $h\in\mathcal E_\sigma^N$ if
and only if it is a primitive function of an
element in $\mathcal E_\sigma^{N-1}$. The
relationship between $\mathcal E_\sigma^N$
and $\mathcal H_{\alpha,p}^N$ is as follows.

\begin{Lem}If $h\in\mathcal
H_{\alpha,p}^N$, then the restriction of $h$
on $\R\setminus (-1,1)$ is in $\mathcal
E_\sigma^{N}$.
\end{Lem}
\begin{proof} It suffices to show that the restriction of $h^{(N)}$ on $\R\setminus (-1,1)$ is in 
$\mathcal E_\sigma$. This is obvious since  $h^{(N)}_c$ has at most polynomial increasing speed at infinity.
\end{proof}

\begin{Def}\label{Def: lambda h}For any derivable function $h$ on $\R\setminus (-1,1)$, define the operator
\begin{equation}\label{Equ:operator lambda}
\Lambda(h) (x):= \Big(\frac{h(x)}{x}\Big)'.\end{equation}
\end{Def}
\begin{Lem}\label{Lem:f_Gamma} 
If $h\in\mathcal E_\sigma^1$, then $\Lambda(h)\in\mathcal E_\sigma$. Furthermore, we have the following equality:
\begin{equation}\label{Equ:derive de fh}
\widetilde f_h'(x)=x\widetilde f_{\Lambda(h)}(x).
\end{equation}
\end{Lem}
\begin{proof}By definition, $\Lambda(h)(x)=h'(x)/x-h(x)/x^2$, so $\Lambda(h)\in\mathcal E_\sigma$.
To prove the equality, it suffices to verify
that the function $u(x):=x^{-1}\widetilde
f_h'(x)$ satisfies the equation
\eqref{Equ:modified Stein equatli} for
$\Lambda(h)$ (see Remark above). In fact, if we divide the
both side of the equation $x\widetilde
f_h(x)-\sigma^2 \widetilde f_h'(x)=h(x)$ by
$x$ and then take the derivative, we obtain
$xu(x)-\sigma^2u'(x)=\Lambda(h)(x)$.
%Since $\widetilde
%f_{\Lambda(g)}$ is the only solution of this equation
%having at most polynomial increasing speed at infinity,
%we have $u=\widetilde f_{\Lambda(g)}$.
\end{proof}

\begin{Lem}\label{Lem:croissance de f_h}
If $h\in\mathcal E_\sigma$ and if $l$ is a real number
such that $h(x)=O(|x|^l)$, then $\widetilde
f_h(x)=O(|x|^{l-1}) $.
\end{Lem}
\begin{proof}Recall that (\cite[ LemA.1]{EJ2008}) if $|h(x)|\leq g(x)$ and if $g(x)/|x|$ is decreasing when $x>0$ and is increasing when $x<0$, then $|\widetilde f_h(x)|\leq g(x)/|x|$. Hence, we prove the lemma for the cases where $l<1$.
By Lemma \ref{Lem:f_Gamma}, one has
\[\begin{split}
\widetilde
f_{|x|^l}(x)&=x^{-1}(|x|^l+\sigma^2\widetilde
f_{|x|^l}'(x)
)=\mathrm{sgn}(x)|x|^{l-1}+\sigma^2\widetilde
f_{\Lambda(|x|^l)}(x)\\
&=\mathrm{sgn}(x)|x|^{l-1}+\sigma^2(l-1)\widetilde
f_{|x|^{l-2}}(x).
\end{split}\]Thus,
$\widetilde f_{|x|^{l-2}}=O(|x|^{l-3})$
implies $\widetilde
f_{|x|^{l}}=O(|x|^{l-1})$. Hence by induction on $l$, we obtain the result.
\end{proof}

\begin{Rem}\label{Rem: bounded tilde fh by polynomial}
With the notation of Barbour \cite{Ba1986}, the equivalent expectation form of $\widetilde f_h$ is given by
\begin{equation}\label{Equ:widetilde f expectation form}
\widetilde f_{h}(x)=
\begin{cases}
\frac{\sqrt{2\pi}}{\sigma}
\esp\big[h(Z+x)e^{-\frac{Zx}{\sigma^2}}\indic_{\{Z>0\}}
\big],\quad
&x>0\\
-\frac{\sqrt{2\pi}}{\sigma}
\esp\big[h(Z+x)e^{-\frac{Zx}{\sigma^2}}\indic_{\{Z<0\}}\big],\quad
&x<0
\end{cases}
\end{equation}
where $Z\sim N(0,\sigma^2)$.
So the above lemma can be interpreted as~: the function
\[\frac{1}{x^l}\mathbb E\big[\indic_{\{Z>0\}}(Z+x)^{l+1}e^{-\frac{Zx}{\sigma^2}}\big]\]
is bounded on $[1,+\infty)$. We can then deduce easily the following assertion~: for all $l\in\R$ and $m\in\R_+$, the function
\[\frac{1}{x^l}\mathbb E\big[\indic_{\{Z>0\}}(Zx)^m(Z+x)^{l+1}e^{-\frac{Zx}{\sigma^2}}\big]\]
is bounded on $[1,+\infty)$ by using the fact that the function $u^me^{-\frac{u}{2\sigma^2}}$ is bounded 
on $[0,\infty)$.
\end{Rem}

%=====================================
\def\skip{
 Since we are interested in the
estimations of $\|\cdot\|_{\alpha,p}^A$ with $A=(-\infty,1]$ or
$A=[1,+\infty)$, such truncation makes no harm. In the following, we
shall define a sequence of subspaces of $\mathcal E_\sigma$ and some
operators acting on them. We say a Borel function $g$ on
$(-\infty,1]\cup[1,+\infty)$ has at most polynomial increasing speed
at infinity if there exists $a\in\mathbb N$ such that the function
$g(x)/x^a$ is bounded. Denote by $\mathcal E_\sigma^0$ the space of
such functions. One has $\mathcal E_\sigma^0\subset\mathcal
E_\sigma$. For any integer $N$, $N\ge 1$, let $\mathcal E_\sigma^N$
be the space of $(N-1)^{\mathrm{th}}$ order differentiable function
$h$ on $(-\infty,1]\cup[1,+\infty)$ such that $h^{(N-1)}$ is
absolutely continuous and is a primitive function of an element in
$\mathcal E_\sigma^0$. In the following is some properties on these
function spaces.

\begin{Pro} \label{Pro:properties of e}Let $N$ be an integer $N\ge 0$.
\begin{enumerate}[1)]
\item For any integer $N\ge 1$, $\mathcal E_\sigma^N\subset
\mathcal E_\sigma^{N-1}$. Moreover, for any absolutely continuous
function $h$, $h\in\mathcal E_\sigma^N$ if and only if it is a
primitive function of an element in $\mathcal E_\sigma^{N-1}$.
\item The space $\mathcal E_\sigma^N$ is stable by
multiplication by Laurent polynomials.
\end{enumerate}
\end{Pro}
\begin{proof}
1) The second assertion results by definition. We shall prove the
first one. Actually, we only need to verify $\mathcal
E_\sigma^1\subset\mathcal E_\sigma^0$. Assume that $h\in\mathcal
E_{\sigma}^1$, the it has at most polynomial increasing speed at
infinity since it is a primitive function of a function of this
type. It suffices then to verify that $h\in\mathcal E_\sigma$.
Denote by $h'\in\mathcal E_\sigma^0$ an density function of $h$ with
respect to the Lebesgue measure. Let $P$ be an arbitrary polynomial.
One has
\[\int_1^{+\infty}|P(x)h(x)|\phi_\sigma(x)\,dx\le
\int_1^{+\infty}dx\,|P(x)|\phi_\sigma(x)\int_1^x|h'(t)|\,
dt+|h(1)|\int_1^{+\infty}|P(x)|\phi_\sigma(x)\,dx.\]
The last integral is finite, and the double integral
equals
\[\int_1^{+\infty}dt\,|h'(t)|\int_t^{+\infty}
|P(x)|\phi_\sigma(x)\,dx\] By integration by
part, there exists a polynomial $Q(x)$ such that
$\int_t^{+\infty}|P(x)|\phi_\sigma(x)\,dx\le
Q(t)$ for any $t\ge 1$, which implies that the double
integral above is finite since $h'\in\mathcal
E^0_\sigma$.

2) Let $h\in\mathcal E_\sigma^N$, $R$ be a Laurent
polynomial, and $g=Rh$. The function $g^{(N)}$ can be
written as linear combination of $h,h',\cdots,h^{(N)}$
with Laurent polynomial coefficients. From 1), we
obtain that these derivatives of $h$ are all in
$\mathcal E_\sigma^0$, hence also is $g^{(N)}$.
\end{proof}
\begin{Def}
If $h$ is an element in $\mathcal E_\sigma^1$, denote
by $\Lambda(h)$ the function on $\mathbb R\setminus(-1,1)$
defined as $\Lambda(h)(x)=(h(x)/x)'$.
\end{Def}

\begin{Pro}
The mapping $\Lambda$ is a linear operator from $\mathcal
E_\sigma^1$ to $\mathcal E_\sigma^0$. More generally,
its restriction on $\mathcal E_\sigma^N$ sends
$\mathcal E_\sigma^N$ to $\mathcal E_\sigma^{N-1}$.
\end{Pro}
\begin{proof}
By definition, $\Lambda(h)(x)=h'(x)/x-h(x)/x^2$. So the
assertion results from Proposition \ref{Pro:properties
of e}.
\end{proof}
}

%=======================================

We give below the relationship between the
derivatives of $\widetilde f_h$ and of $h$.
In the following two formulas, the first one
computes $\widetilde f_h^{(N)}$ using the
operator \eqref{Equ:operator lambda} and the
second one expresses $\Lambda^{N}(h)$ using
derivatives of $h$. Their proofs are by
induction, which we omit in this
article (interest readers may refer to
\cite[p.144-145]{Ji2006}). We only remind that the first
formula is a generalization  of
\eqref{Equ:derive de fh}.

\begin{Lem}\label{Pro:f_h^{(n)} as gamma^n(h), normal}
If $h\in\mathcal E^N_\sigma$ with $N$ being a strictly positive integer,  then
\begin{gather}
\label{Equ:higher derivatives of fh} \widetilde
f_h^{(N)}(x)=\sum_{k=0}^{\lfloor N/2\rfloor}\binom{N}{2k}
(2k-1)!!x^{N-2k}\widetilde f_{\Lambda^{N-k}(h)}(x);\\
\label{Equ:Gamma N(h)} \Lambda^{N}(h)(x)=\sum_{k=0}^{N}
(-1)^{k}
(2k-1)!!\binom{N+k}{2k}\frac{h^{(N-k)}(x)}{x^{N+k}}.
\end{gather}
where we have used the convention $(-1)!!=1$
and  $\lfloor N/2\rfloor$ denotes the largest
integer not exceeding ${N}/{2}$.
\end{Lem}

\begin{Rem}\label{Rem:combinaison linearire}
\begin{enumerate}[1)]
\item For any function $h\in\mathcal E_\sigma^N$, the above results also hold for $\widetilde
f_h^{(m)}(x)$ and $\Lambda^m(h)$ where $1\leq m\leq N$. As
the operator $h\rightarrow
\widetilde f_h$ is linear on $h$,  the
above lemma enables us to write the
derivatives of $\widetilde f_h$ as a linear
combination of derivatives of $h$ with
Laurent polynomial coefficients and then to  deduce their
increasing speed at infinity.

\item The derivative function $\widetilde f^{(N+1)}_h$ has to be treated differently.
In fact, we can no longer apply \eqref{Equ:Gamma N(h)} to $\Lambda^{N+1}(h)$ since $h^{(N+1)}$ does 
not necessarily exist. We separate the first term where $k=0$ in \eqref{Equ:higher derivatives of fh} from the others and then take the derivative to obtain
\[\begin{split}\widetilde f^{(N+1)}_h&=x^N\widetilde{f}'_{\Lambda^N(h)}(x)+\sum_{k=0}^{\lfloor N/2\rfloor}\binom{N}{2k}
(2k-1)!!(N-2k)x^{N-2k-1}\widetilde f_{\Lambda^{N-k}(h)}(x)\\&\qquad\qquad\quad\quad+\sum_{k=1}^{\lfloor N/2\rfloor}\binom{N}{2k}
(2k-1)!!x^{N-2k}\widetilde f'_{\Lambda^{N-k}(h)}(x)\\
&=x^N\widetilde{f}'_{\Lambda^N(h)}(x)+\sum_{k=1}^{\lfloor (N+1)/2\rfloor}\binom{N+1}{2k}
(2k-1)!!x^{N+1-2k}\widetilde f_{\Lambda^{N+1-k}(h)}(x).
\end{split}\]
This will be a crucial point in the proof of Proposition
\ref{Pro:finitude}.
\end{enumerate}
\end{Rem}

\begin{Lem}\label{Lem:B14}
Let $h$ be a Borel function defined on
$A_2\cup A_3=\mathbb R\setminus (-1,1)$. If
$\|h\|_{\alpha,p}^A<+\infty$ where $A=A_2$ or $A_3$,
then, for any integer $n\geqslant 0$, one has
\[\|x^{n+1}\widetilde
f_{{h}/{x^{n}}}\|_{\alpha,p}^A<+\infty,\qquad
\|x^n\widetilde
f_{{h}/{x^n}}'\|_{\alpha,p}^A<+\infty .\]
\end{Lem}
\begin{proof}We only prove for $A_2$ and the case for $A_3$ is by symmetry. 
Let $g(x)=h(x)/x^n$. Assume that $x$ and $y$ are two
real numbers such that $1\leqslant x<y$. Then one has
\[\frac{|x^{n+1}\widetilde f_g(x)-y^{n+1}\widetilde{f}_g(y)|}{
|x-y|^\alpha(1+|x|^p+|y|^p)}=
\frac{\sqrt{2\pi}}{\sigma}\mathbb E\bigg[
I_{\{Z>0\}}\frac{\Big|\frac{h(Z+x)x^{n+1}}{(Z+x)^{n}}e^{-{Zx}/{\sigma^2}}
-\frac{h(Z+y)y^{n+1}}{(Z+y)^{n}}e^{-{Zy}/{\sigma^2}}\Big|}{|x-y|^\alpha(1+|x|^p+|y|^p)}
\bigg],\] which can be bounded from above by the sum of the
following two terms
\begin{gather}\label{Equ:formule 23}\frac{\sqrt{2\pi}}{\sigma}\mathbb
E\bigg[I_{\{Z>0\}}\frac{|h(Z+x)-h(Z+y)|}{|x-y|^\alpha(1+|x|^p+|y|^p)}
\cdot\frac{y^{n+1}}{(Z+y)^{n}}e^{-\frac{Zy}{\sigma^2}}\bigg]\\
\label{Equ:formule 24}\frac{\sqrt{2\pi}}{\sigma}\mathbb
E\bigg[I_{\{Z>0\}}|h(Z+x)|
\frac{\Big|\frac{x^{n+1}}{(Z+x)^{n}}e^{-\frac{Zx}{\sigma^2}}
-\frac{y^{n+1}}{(Z+y)^{n}}e^{-\frac{Zy}{\sigma^2}}\Big|}{
|x-y|^\alpha(1+|x|^p+|y|^p)}\bigg]
\end{gather}
Note that \eqref{Equ:formule 23} is bounded from above
by
\[\|h\|_{\alpha,p}^{A_2}\frac{\sqrt{2\pi}}{\sigma}
y^{n+1}\mathbb
E\Big[I_{\{Z>0\}}\frac{1}{(Z+y)^{n}}e^{-\frac{Zy}{\sigma^2}}\Big]
=\|h\|_{\alpha,p}^{A_2}\,y^{n+1}\widetilde
f_{\frac 1{|x|^n}}(y).\] By Lemma
\ref{Lem:croissance de f_h}, this quantity is bounded.
We then consider the upper bound of
\eqref{Equ:formule 24} under the supplementary
condition that $y\leqslant 2x$. As
$\|h\|_{\alpha,p}^{A_2}<+\infty$, there exists
a constant $C>0$ such that $h(x)\leqslant
C|x|^{\alpha+p}$. By applying the mean value theorem on the function $\frac{x^{n+1}}{(Z+x)^n}e^{-\frac{Zx}{\sigma^2}}$ and the fact that
$|x-y|^{\alpha-1}(1+|x|^p+|y|^p)\geq |x|^{\alpha-1+p}$ where $\alpha<1$, 
\[\begin{split}\eqref{Equ:formule 24}&\leqslant \frac{\sqrt{2\pi}}{\sigma}
|x|^{1-\alpha-p}\mathbb E\bigg[I_{\{Z>0\}}C(Z+x)^{\alpha+p}
e^{-\frac{Zx}{\sigma^2}}\Big(\frac{(n+1)(2x)^{n}}{(Z+x)^{n}}
+\frac{n(2x)^{n+1}}{(Z+x)^{n+1}}+\frac{Z(2x)^{n+1}}{\sigma^2(Z+x)^{n}}\Big)\bigg]\\
&\leqslant C
\frac{\sqrt{2\pi}}{\sigma}|x|^{1-\alpha-p}\bigg\{\big(2^n\cdot(n+1)+2^{n+1}\cdot
n\big)\mathbb
E\big[I_{\{Z>0\}}(Z+x)^{\alpha+p}e^{-\frac{Zx}{\sigma^2}}\big]\\&\qquad\qquad\qquad\qquad+2^{n+1}
\mathbb
E\Big[I_{\{Z>0\}}\frac{Zx}{\sigma^2}(Z+x)^{\alpha+p}e^{-\frac{Zx}{\sigma^2}}\Big]\bigg\}\\
&\ll
|x|^{1-\alpha-p}\Big\{E\big[I_{\{Z>0\}}(Z+x)^{\alpha+p}e^{-\frac{Zx}{\sigma^2}}\big]+
\mathbb
E\Big[I_{\{Z>0\}}\frac{Zx}{\sigma^2}(Z+x)^{\alpha+p}e^{-\frac{Zx}{\sigma^2}}\Big]
\Big\},
\end{split}
\]
which is bounded (see Remark \ref{Rem: bounded tilde fh by polynomial}). In the case where $y>2x$, one has $x^{n+1}/(Z+x)^n\leq x$ when $Z\geq 0$, and $|x-y|^{\alpha}(1+|x|^p+|y|^p)\geq\max(|x|^{\alpha+p},\big(\frac{y}{2}\big)^\alpha\cdot y^p)$, so
\[\eqref{Equ:formule 24}\leqslant C\frac{\sqrt{2\pi}}{\sigma}
\left\{ \frac{1}{2|x|^{\alpha+p-1}}\mathbb E
\Big[I_{\{Z>0\}}
(Z+x)^{\alpha+p}e^{-\frac{Zx}{\sigma^2}}\Big]+
\frac{2^\alpha}{|y|^{\alpha+p-1}}\mathbb E
\Big[I_{\{Z>0\}}
(Z+y)^{\alpha+p}e^{-\frac{Zy}{\sigma^2}}\Big]\right\},\]
which is bounded.
\end{proof}

We now give the final part of the proof.
\begin{Pro}\label{Pro:finitude}
Let $N\geqslant 0$ be an integer,
$\alpha\in(0,1]$ and $p\geqslant 0$. Let $h$ be a function defined on $A_2\cup A_3=\mathbb
R\setminus(-1,1)$ which is
$N$ times differentiable and such that
$h^{(N)}$ is locally of finite variation,
having finitely many jumps and verifying
$\|h^{(N)}_c\|_{\alpha,p}^{A}<+\infty$ where $A=A_2$ or $A_3$.
Then the function $\widetilde f_h$ is
$N+1$ times differentiable,
$\widetilde f_h^{(N+1)}$ is locally of finite
variation, having finitely many jumps
and verifying
$\|\widetilde
f_{h,c}^{(N+1)}\|_{\alpha,p}^{A}<+\infty.$
\end{Pro}
\begin{proof}
The function $\widetilde f_h$ is
$N+1$ times differentiable by \eqref{Equ:derive d'ordre N de Stein equation}, $\widetilde
f_h^{(N+1)}$ is locally of finite variation, having
only finitely many jumps. By virtue of
Lemmas \ref{Lem:two intervals} and \ref{Pro:bounded on
finite interfval}, it suffices to prove
\[\max\{\|\widetilde
f_{h,c}^{(N+1)}\|_{\alpha,p}^{(-\infty,-b]},
\|\widetilde
f^{(N+1)}_{h,c}\|_{\alpha,p}^{[b,+\infty)}\}<+\infty\]
for sufficiently positive number $b$.
Therefore, without loss of generality, we may
assume that $h^{(N)}$ is continuous and hence
$\widetilde f_h^{(N+1)}$ is also continuous.\\
By Remark \ref{Rem:combinaison linearire} 2), the function
$\widetilde f^{(N+1)}_h$ can be written as a linear
combination of  $x^N\widetilde f'_{\Lambda^N(h)}$ and
terms of the form $x^{N+1-2k}\widetilde f_{\Lambda^{N+1-k}(h)}(x)$ where $k= 1,\cdots,\lfloor\frac{N+1}{2}\rfloor$. By \eqref{Equ:Gamma N(h)}, $x^N\widetilde f'_{\Lambda^N(h)}$ itself is also a linear combination of $x^N\widetilde f'_{h^{(N-i)}/x^{N+i}}$ where $i=0,\cdots, N$. 
As $\|h^{(N)}\|^A_{\alpha,p}<\infty$,  we have, similar as in in Lemma \ref{Lem:Lipschitz norms} 4), that $\|h^{(N-i)}\|_{\alpha,p+i}^{A}<+\infty$. Hence $\|h^{(N-i)}/x^{i}\|^A_{\alpha,p}<\infty$ by Lemma \ref{Lem:5.2} and Lemma \ref{Lem:B14}  then implies that $\|x^N\widetilde f'_{h^{(N-i)}/x^{N+i}}\|^A_{\alpha,p}<\infty$.\\
The terms $x^{N+1-2k}\widetilde f_{\Lambda^{N+1-k}(h)}(x)$ are also, by \eqref{Equ:Gamma N(h)} again, linear combinations of the functions of the form $x^{N+1-2k}\widetilde f_{{h^{(N+1-k-i)}}/{x^{N+1-k+i}}}$.  By a similar argument as above using Lemma  \ref{Lem:5.2}, $\|h^{(N+1-k-i)}/x^{1+k+i}\|_{\alpha,p}^A<\infty$. Finally, we apply Lemma \ref{Lem:B14} to obtain $\|x^{N+1-2k}\widetilde f_{{h^{(N+1-k-i)}}/{x^{N+1-k+i}}}\|^A_{\alpha,p}<\infty$, which completes the proof.

\end{proof}

\bibliography{jiao}

\begin{thebibliography}{10}

\bibitem{Ba1986}
A.~D. Barbour.
\newblock Asymptotic expansions based on smooth functions in the central limit
  theorem.
\newblock {\em Probability Theory and Related Fields}, 72:289--303, 1986.

\bibitem{Ba1987}
A.~D. Barbour.
\newblock Asymptotic expansions in the poisson limit theorem.
\newblock {\em Annals of Probability}, 15:748--766, 1987.

\bibitem{CS2001}
L.~H.~Y. Chen and Q.-M. Shao.
\newblock A non-uniform {B}erry-{E}sseen bound via {S}tein's method.
\newblock {\em Probability Theory and Related Fields}, 120:236--254, 2001.

\bibitem{CS2005}
L.~H.~Y. Chen and Q.-M. Shao.
\newblock Stein's method for normal approximation.
\newblock In {\em An Introduction to {S}tein's Method}, volume~4 of {\em
  Lecture Notes Series, IMS, National University of Singapore}, pages 1--59.
  Singapore University Press and World Scientific Publishing Co. Pte. Ltd.,
  2005.

\bibitem{EJ2008}
N.~El~Karoui and Y.~Jiao.
\newblock Stein's method and zero bias transformation for {C}{D}{O}s tranches
  pricing.
\newblock to appear in {F}inance and {S}tochastics, 2008.

\bibitem{Go2007}
L.~Goldstein.
\newblock ${L}^1$ bounds in normal approximation.
\newblock {\em Annals of Probability}, 35(5):1888--1930, 2007.

\bibitem{GR1997}
L.~Goldstein and G.~Reinert.
\newblock {S}tein's method and the zero bias transformation with application to
  simple random sampling.
\newblock {\em Annals of Applied Probability}, 7:935--952, 1997.

\bibitem{GH1978}
F.~G{\"o}tze and C.~Hipp.
\newblock Asymptotic expansions in the central limit theorem under moment
  conditions.
\newblock {\em Z. Wahrscheinlichkeitstheorie und Verw. Gebiete}, 42(1):67--87,
  1978.

\bibitem{Hi1977}
C.~Hipp.
\newblock Edgeworth expansions for integrals of smooth functions.
\newblock {\em Ann. Probability}, 5(6):1004--1011, 1977.

\bibitem{Ji2006}
Y.~Jiao.
\newblock Risque de cr\'edit: mod\'elisation et simulation num\'erique.
\newblock PhD thesis, Ecole Polytechnique,
  http://www.imprimerie.polytechnique.fr/Theses/Files/Ying.pdf, 2006.

\bibitem{KF1974}
A.~Kolmogolov and S.~Fomine.
\newblock {\em \'El\'ements de la th\'eorie des fonctions et de l'analyse
  fonctionnelle}.
\newblock \'Editions Mir., Moscow, 1974.
\newblock Avec un compl\'ement sur les alg\`ebres de Banach, par V. M.
  Tikhomirov, Translated from Russia to French by Michel Dragnev.

\bibitem{La2001}
F.~Laudenbach.
\newblock {\em Calcul diff\'erentiel et int\'egral}.
\newblock \'Editions \'Ecole Polytechnique, 2001.

\bibitem{Pe1975}
V.~V. Petrov.
\newblock {\em Sums of Independent Random Variables}.
\newblock Springer-Verlag, 1975.

\bibitem{Ro2005}
V.~Rotar.
\newblock Stein's method, {E}dgeworth's expansions and a formula of {B}arbour.
\newblock In {\em Stein's method and applications}, volume~5 of {\em Lect.
  Notes Ser. Inst. Math. Sci. Natl. Univ. Singap.}, pages 59--84. Singapore
  Univ. Press, Singapore, 2005.

\bibitem{St1972}
C.~Stein.
\newblock A bound for the error in the normal approximation to the distribution
  of a sum of dependent random variables.
\newblock In {\em Proc, Sixth Berkeley Symp. Math. Statist. Probab.}, pages
  583--602. Univ. California Press, Berkeley, 1972.

\bibitem{St1986}
C.~Stein.
\newblock {\em Approximate Computation of Expectations}.
\newblock IMS, Hayward, CA., 1986.

\end{thebibliography}
\bibliographystyle{plain}

\end{document}